\documentclass[12pt,reqno]{amsart}

\usepackage[activeacute,english]{babel}
\usepackage[utf8]{inputenc}
\usepackage{amsfonts,amsmath,amsthm,bm,amssymb}
\usepackage{graphics,tikz,float}

\usepackage{hyperref,cleveref}
\usepackage{autonum}

\DeclareRobustCommand{\SkipTocEntry}[5]{}

\usepackage[margin=0.9in]{geometry}
\newcommand{\BB}{{\mathcal B}}

\newcommand{\PP}{{\mathcal P}}

\newcommand{\R}{{\mathbb R}}
\newcommand{\Z}{{\mathbb Z}}

\newcommand{\Sph}{{\mathbb S}}

\newcommand{\eper}{\operatorname{per}}

\newtheorem{theorem}{Theorem}[section]
\newtheorem{lemma}[theorem]{Lemma}

\newtheorem{remark}[theorem]{Remark}

\numberwithin{equation}{section}

\allowdisplaybreaks[1]

\subjclass[2010]{47G20,35J60,35B10,49Q05}
\keywords{Constant nonlocal mean curvature, periodic surfaces, Delaunay cylinders}

\begin{document}

\title[Periodic nonlocal-CMC surfaces]
{Existence and symmetry of periodic nonlocal-CMC \\ surfaces
via variational methods}

\author[X. Cabré]{Xavier Cabr\'e}
\address{X. Cabr\'e \textsuperscript{1,2,3}
\newline
\textsuperscript{1} ICREA, Pg. Lluis Companys 23, 08010 Barcelona, Spain
\newline
\textsuperscript{2} Universitat Polit\`ecnica de Catalunya, Departament de Matem\`{a}tiques and IMTech, Diagonal 647, 08028 Barcelona, Spain
\newline
\textsuperscript{3}
Centre de Recerca Matem\`atica, Edifici C, Campus Bellaterra, 08193 Bellaterra, Spain.}
\email{xavier.cabre@upc.edu}

\author[G. Csat\'o]{Gyula Csat\'o}
\address{G. Csat\'o \textsuperscript{1,2}
\newline
\textsuperscript{1}
Departament de Matem\`{a}tiques i Inform\`{a}tica,
Universitat de Barcelona,
Gran Via 585,
08007 Barcelona, Spain
\newline
\textsuperscript{2}
Centre de Recerca Matem\`atica, Edifici C, Campus Bellaterra, 08193 Bellaterra, Spain.}
\email{gyula.csato@ub.edu}

\author[A. Mas]{Albert Mas}
\address{A. Mas \textsuperscript{1,2}
\newline
\textsuperscript{1}
Departament de Matem\`atiques,
Universitat Polit\`ecnica de Catalunya,
Campus Diagonal Bes\`os, Edifici A (EEBE), Av. Eduard Maristany 16, 08019
Barcelona, Spain
\newline
\textsuperscript{2}
Centre de Recerca Matem\`atica, Edifici C, Campus Bellaterra, 08193 Bellaterra, Spain.}
\email{albert.mas.blesa@upc.edu}

\date{\today}
\thanks{The three authors are supported by the Spanish grants MTM2017-84214-C2-1-P, 
PID2021-123903NB-I00, and RED2018-102650-T funded by MCIN/AEI/10.13039/501100011033 and by ``ERDF A way of making Europe''. The second and third authors are also supported by the Spanish grant MTM2017-83499-P. The second author is in addition supported by the Spanish grant PID2021-125021NA-I00.
This work is supported by the Spanish State Research Agency, through the Severo Ochoa and Mar\'ia de Maeztu Program for Centers and Units of
Excellence in R\&D (CEX2020-001084-M)}

\begin{abstract}
This paper provides the first variational proof of the existence of periodic nonlocal-CMC surfaces. These are nonlocal analogues of the classical Delaunay cylinders. 
More precisely, we show the existence of a set in $\R^n$
which is periodic in one direction, has a prescribed (but arbitrary) volume within a slab orthogonal to that direction, has constant nonlocal mean curvature, and minimizes an appropriate periodic version of the fractional perimeter functional under the volume constraint. We show, in addition, that the set is cylindrically symmetric and, more significantly, that it is even as well as nonincreasing on half its period. This monotonicity property solves an open problem and an obstruction which arose in an earlier attempt, by other authors, to show the existence of minimizers.
 \end{abstract}

\maketitle



\section{Introduction}

For a general set $E\subset\R^n$, the nonlocal (or $s$-fractional, $0<s<1$) mean curvature of
$\partial E$ (or of $E$) at a boundary point $x\in\partial E$ is defined by
\begin{equation}\label{NMC}
H_{s}[E](x) := PV \int_{\R^n}dy\,\frac{\chi_{E^c}(y)-\chi_{E}(y)}{|x-y|^{n+s}}
=\lim_{\epsilon\downarrow 0}\int_{\{|x-y|>\epsilon\}´}dy\,\frac{\chi_{E^c}(y)-\chi_{E}(y)}{|x-y|^{n+s}}.
\end{equation}
Here $PV$ denotes the principal value sense, $\chi_E$ is the characteristic function of the set $E,$ and $E^c$ is the complement in $\R^n$ of the set $E$.  The nonlocal mean curvature arises when computing the first variation of the fractional perimeter, an extension of the classical perimeter that was introduced, in its localized or Dirichlet version, in the pioneering article of Caffarelli, Roquejoffre, and Savin \cite{CafRoqSav}. Their work has allowed for a natural extension of the theory of minimal surfaces to the nonlocal framework.

The current paper concerns the nonlocal analogue of the classical result of Delaunay \cite{Delaunay 1841} on the existence of periodic surfaces of revolution with constant mean curvature (CMC). 
We provide the first variational proof of the existence of periodic surfaces of revolution having constant nonlocal mean curvature  and a prescribed (but arbitrary) volume within one period. Throughout the paper, surfaces with constant nonlocal mean curvature will be called nonlocal-CMC surfaces.
Before our work, and up to our knowledge, the only papers which had established the existence of some periodic nonlocal-CMC surfaces are \cite{CFSW,CFW1,CFW2}, by Cabr\'e, Fall, Sol\`a-Morales, and Weth. These works
did not use variational methods,
but perturbative techniques (the implicit function theorem, essentially). Thus, they only obtained nonlocal-CMC surfaces which are very close to certain explicit configurations, as explained next.  

In arbitrary dimension, \cite{CFW1} showed the existence of a one parameter family of periodic nonlocal-CMC surfaces bifurcating from a straight cylinder. The paper extended a previous analogous result in $\R^2$ from~\cite{CFSW}.\footnote{Such nonlocal-CMC curves in the plane have no analogue, obviously, for the classical mean curvature.} Later, also by perturbation techniques,
 \cite{CFW2} proved the existence of a periodic array of disjoint near-balls, the full array having constant nonlocal mean curvature. They bifurcate from infinity, starting from a single round ball of radius $1$ ---in particular, the near-balls have all the same shape (which is very close to a unit ball) and are very distant from each other. Thus, after rescaling them to have their consecutive centers at distance $2\pi$, they become of very small size (and, hence, of small volume within a period). Instead, the $2\pi$-periodic near-cylinders from \cite{CFW1} have a certain positive, bounded below, width. Prior to the current paper, no result connecting both configurations, through other periodic nonlocal-CMC surfaces, was available.\footnote{Recall that, in the local case, the classical Delaunay surfaces vary continuously from a cylinder to an infinite compound of tangent spheres.}

This paper establishes for the first time the existence of a periodic nonlocal-CMC surface in~$\R^n$ for every given volume within a slab. 
By a variational method 
we show the existence of a minimizer of the functional
\begin{equation}\label{p} 
\PP_s[E] := \int_{E\cap\{-\pi<x_1<\pi\}}\!dx
\int_{\R^n\setminus E}dy\, \frac{1}{|x-y|^{n+s}}
\end{equation}
(introduced for the first time in \cite{Davila}),
under a volume constraint within the slab 
$\{x\in\R^n:\,-\pi<x_1<\pi\}$, among subsets $E$ of $\R^{n}=\{x=(x_1,x')\in\R\times\R^{n-1}\}$ which are $2\pi$-periodic in the coordinate $x_1$. 
Notice that the functional $\PP_s$ differs from the fractional $s$-perimeter of $E$ relative to 
$\{x\in\R^n:\,-\pi<x_1<\pi\}$, a quantity which is given by 
adding an additional term to $\PP_s$, as follows: 
\begin{equation}\label{DirPer}
  \PP_s[E]  + \int_{(\R^n\setminus E)\cap\{-\pi<x_1<\pi\}}\!dx
  \int_{E\cap\{|y_1|>\pi\}}dy\, \frac{1}{|x-y|^{n+s}}.
\end{equation}
Still, and somehow surprisingly, 
we establish that the first variation of these two functionals ---the periodic first variation for \eqref{p} and the Dirichlet first variation in the slab for \eqref{DirPer}--- lead both to the nonlocal mean curvature.\footnote{The Dirichlet first variation of \eqref{DirPer} was derived in the seminal paper \cite{CafRoqSav} by Caffarelli, Roquejoffre, and Savin.} Uncovering this fact originated the current article.

The following is our first main result. Apart from the existence of a minimizer, it also establishes symmetry and monotonicity properties of any minimizer.

\begin{theorem}\label{scan thm1}
Given $n\geq2$, $s\in(0,1)$, and $\mu>0$, there exists a minimizer $E\subset\R^n$ of $\PP_s$ among all measurable sets $F\subset\R^n$ which are $2\pi$-periodic in the $x_1$-variable and satisfy
\begin{equation}
  |F\cap\{-\pi<x_1<\pi\}|=\mu.
\end{equation}
In addition, up to sets of measure zero:
\begin{itemize}
\item[$(i)$] Up to a translation in the 
$x'$-variable, every minimizer $E$ is of the form
\begin{equation}\label{Del frac min E u}
E=\{x=(x_1,x')\in\R\times\R^{n-1} :\,|x'|<u(x_1)\}
\end{equation}
for some $2\pi$-periodic nonnegative function $u$ whose power $u^{n-1}$ belongs to $W^{s,1}(-\pi,\pi)$. In particular, we have
$u\in W^{\frac{s}{n-1},\,n-1}(-\pi,\pi)$.

\item[$(ii)$] Up to a translation in the 
$x_1$-direction, $u$ is even, as well as nonincreasing in $(0,\pi)$. 

\item[$(iii)$] For every minimizer $E$, there exists a constant $\lambda\in\R$ such that if $x\in\partial E$, $u$ is $C^{1,\alpha}$ in a neighborhood of $x_1$ for some $\alpha>s$,  and $u(x_1)>0$, then $H_s[E](x)=\lambda$.

\item[$(iv)$] For $\mu$ small enough depending only on $n$ (and hence uniformly as $s\uparrow1$ and as $s\downarrow0$), every minimizer is not a straight circular cylinder, i.e., the function $u$ is nonconstant.
\end{itemize}

\end{theorem} 

The symmetry and monotonicity properties of points $(i)$ and $(ii)$ will be discussed in the next Subsection~\ref{symmetry}. Of particular importance is point $(ii)$, since it relies on a deep and not so well-known result on rearrangements in spheres $\Sph^{m}$. However, we emphasize that our proof of existence of a minimizer only uses simple tools and, in particular, it is independent of this rearrangement result.

The work \cite{Davila} was the first attempt to show the existence of periodic nonlocal-CMC surfaces with a minimization technique. This paper introduced, for the first time, the functional $\PP_s$ in \eqref{p} for periodic sets. However, it did not relate its Euler-Lagrange equation to the nonlocal mean curvature.\footnote{At that time, we knew that the first variation of $\PP_s$ for periodic sets is the nonlocal mean curvature, since we were working on the analogous functional for periodic solutions to fractional semilinear equations (work of us still to appear). When \cite{Davila} became available, we communicated this fact to its authors, we pointed it out at the beginning of page 35 in the printed version of \cite{CFW1}, and was also mentioned in the Master's Thesis of M. Alviny\`a \cite{Alvinya} directed within our group.} At the same time, \cite{Davila} did not prove the existence of a periodic constrained minimizer of $\PP_s$ ---it only obtained a minimizer within the class of sets which are even in the $x_1$-variable and nonincreasing in $(0,\pi)$. Note that this restricted minimization class (made of nonincreasing functions) prevents from checking the Euler-Lagrange equation at a minimizer.
The result on rearrangements in spheres $\Sph^{m}$, mentioned above in connection with point $(ii)$ of our theorem, is the additional tool that would have completed the proof in \cite{Davila}. Instead, as explained later, \cite{Davila} proved a much stronger result than our statement $(iv)$ (on small volumes) in Theorem~\ref{scan thm1}.

The work \cite{Davila} made apparent, thus, that showing the existence of a periodic constrained minimizer of~$\PP_s$ was a delicate issue. In this article we succeed to find a proof which does not rely, in contrast with the approach in \cite{Davila}, on the nonincreasing property of the functions~$u$ in a minimizing sequence. 
To accomplish this, the main difficulty was finding
a  suitable, more treatable expression for the periodic nonlocal perimeter. 
This expression is given in Lemma~\ref{scan lemma 2}, a result that also contains the important lower bound~\eqref{p4}. With these tools at hand, the existence of minimizer will follow readily from the compact embedding of $W^{s,1}(-\pi,\pi)$ into $L^1(-\pi,\pi)$.

In the course of this work, we first completed our existence proof in the case of the plane, $n=2$. This was a simpler task, since it relied on an expression for the fractional perimeter which is similar to some for the nonlocal mean curvature contained in our 2D results from \cite{CFSW}. Since the expression for $n=2$ is more explicit than the one in Lemma~\ref{scan lemma 2}, we will present it in Appendix~\ref{ss CFSW n2}.

Our proof of point $(iv)$ in Theorem \ref{scan thm1} ---stating that minimizers for small volumes are not straight cylinders, uniformly in $s$--- is simple. It is based on an energy-comparison argument and does not require the knowledge of minimizers being nonincreasing in $(0,\pi)$. Instead, Theorem 3 of the above mentioned paper \cite{Davila} proved a much stronger, quantitative version of this result. Using the monotonicity of the minimizer, it established that for small volume constraints a minimizer must be close, in a measure theoretical sense, to a periodic array of round balls. Notice here that the function $u$ defining the minimizer could be identically zero in a subinterval $(a,\pi)$ of $(0,\pi)$ and that this allows for the minimizer to be a periodic array of near-balls. 
At the same time, the near-balls found in \cite{CFW2} are very close to round balls since they were found through the implicit function theorem.
Thus, the result of \cite{Davila} strongly suggests that for small volumes the  minimizers should agree with the near-balls constructed in \cite{CFW2}. However, this is still an open problem.

Regarding point $(iii)$ in \Cref{scan thm1}, we recall that the nonlocal mean curvature is well defined at points $x\in \partial E$ where $\partial E$ is of class $C^{1,\alpha}$,  for some $\alpha>s$, in a neighborhood of~$x$.
However, a complete $C^{1,\alpha}$ regularity theory for volume constrained minimizers is not yet available.\footnote{In the nonperiodic case, the best result (up to our knowledge) on the regularity of volume constrained minimizers of a certain fractional perimeter functional is stated in \cite{Malchiodi Novaga Pagliardini} (which treats sets of revolution in a half-space), claiming their smoothness except for a countable set of points in the rotational axis. Notice that  the Euler-Lagrange equation of the functional in \cite{Malchiodi Novaga Pagliardini} is closely related, but not equal to, the nonlocal mean curvature.}

The following subsection addresses the symmetry and monotonicity properties of minimizers.

\subsection{Symmetry and monotonicity properties of minimizers}\label{symmetry}

Theorem \ref{scan thm1} contains two results on the symmetry and monotonicity properties of minimizing sets: their cylindrical symmetry, \eqref{Del frac min E u}, and the fact that $u$ is an even function which is nonincreasing in $(0,\pi).$ These properties will follow from two rearrangement inequalities.

The first one concerns the {\em cylindrical rearrangement}, defined as follows. For a measurable set $E\subset\R^n$, consider its sections $E_{x_1}=\{x'\in\R^{n-1}:\, (x_1,x')\in E\}$, where $x_1\in\R$. Define $E^{*\!\operatorname{cyl}}$ as the set whose sections $(E^{*\!\operatorname{cyl}})_{x_1}$ for every $x_1\in\R$ are concentric open balls of $\R^{n-1}$ centered at $0$ and of the same $(n-1)$-dimensional Lebesgue measure as $E_{x_1}$. The fact that the cylindrical rearrangement does not increase the functional $\PP_s$ follows from the Riesz rearrangement inequality, and it was already used in \cite[Proposition 13]{Davila}. It will lead to the existence of a cylindrical symmetric minimizer. A stronger fact is that 
every minimizer is cylindrically symmetric, as stated in \Cref{scan thm1}. This will follow from a strict version of the Riesz rearrangement inequality which analyzes the case of equality; see Lemma \ref{gy:lemma:cylindrical rearr}.  

The second rearrangement that we need to consider is the {\em symmetric decreasing periodic rearrangement}, which is defined in the following way. First, for a measurable set $A\subset\R^n$ consider
$$
  A^*:=\text{Steiner symmetrization of $A$  with respect to }  \{x_1=0\}.
$$
This means that $A^*$ is symmetric with respect to the hyperplane $\{x_1=0\}$ and that $A^*\cap\{x'=c\}$ is an  open interval of the same $1$-dimensional Lebesgue measure as $A\cap\{x'=c\}$ for every $c\in\R^{n-1}.$ Now, given a set $E\subset\R^n$ which is $2\pi$-periodic in the $x_1$-variable we define its symmetric decreasing periodic rearrangement by
$$
  E^{*\!\eper}:=\text{$2\pi$-periodic extension of }\left(E\cap \{-\pi<x_1<\pi\}\right)^*,
$$
where the periodic extension is meant in the coordinate $x_1$. 
It is called ``symmetric decreasing'' since the characteristic function $\chi_{E^{*\!\eper}}$ is nonincreasing in $x_1$ within half of the period, $(0,\pi)$, for each given $x'\in\R^{n-1}$. Note also that if $E$ is cylindrically symmetric, i.e., as in \eqref{Del frac min E u}, then
\begin{equation}
  E^{*\!\eper}=\{x\in\R^n\,:\,|x'|<u^{*\!\eper}(x_1)\},
\end{equation}
where $u^{*\!\eper}$ is the $2\pi$-periodic function which is even, nonincreasing in $(0,\pi)$, and equimeasurable with $u$ (see \Cref{section:rearrangement} for the precise definition).

Our second main result establishes that the periodic fractional perimeter $\PP_s$ in \eqref{p} does not increase under  symmetric decreasing periodic rearrangement. This was left as an open problem in the last paragraph of \cite{Davila}. Our result also characterizes the case of equality.

\begin{theorem}
\label{thm:intro:symmetry of Delaunay}
Let $E\subset\R^n$ be a measurable set which is $2\pi$-periodic in the coordinate $x_1$ and satisfies $\PP_s[E]<+\infty.$ 

Then, $\PP_s[E^{*\!\eper}]\leq \PP_s[E].$ Moreover, if 
$\PP_s[E^{*\!\eper}]= \PP_s[E]$ then, up to a set of measure zero, $E=E^{*\!\eper}+ce_1$ for some $c\in\R$, where $e_1:=(1,0,\ldots,0)\in \R^n$.
\end{theorem}

A key ingredient in the proof of \Cref{thm:intro:symmetry of Delaunay} is a Riesz rearrangement inequality on the circle~$\Sph^{1}$, a result which is not as well-known as the classical Riesz rearrangement inequality on~$\R^n$. Such periodic rearrangement inequality  is originally due to Baernstein and Taylor \cite{Baernstein Taylor} and to Friedberg and Luttinger \cite{FriedLutt}, who proved independently different versions of the inequality around the same time; see more details in the comments following Theorem \ref{thm:gy:sharp Friedberg Luttinger}. The case of equality in the Riesz rearrangement inequality on the circle (and more generally on the sphere~$\mathbb{S}^n$) was first treated by Burchard and Hajaiej in \cite{Burchard Hajaiej}.

\subsection{On the shape of minimizers: stability of straight cylinders}\label{stability}

Let us first recall what is known about the shape of minimizers in the local case, that is, for the classical isoperimetric problem relative to a slab.\footnote{In the local case, the isoperimetric problem in a slab is commonly formulated without a periodicity condition. To relate this setting with the periodic one, one considers an even reflection of the slab minimizer. In this way, a period corresponds to twice the width of the slab.}
It turns out that, for $n\leq 8$, the minimizers are either spheres or straight cylinders, depending on whether the volume constraint is smaller or greater, respectively,  than a critical volume; see \cite{Athanassenas,Vogel} or \cite[Theorem~6]{Hauswirth line group} for $n=3$, and \cite{Pedrosa-Ritore} for $n\leq 8$. 
In particular, the classical Delaunay unduloids are not minimizers for $n\leq8$. Furthermore, they are known to be unstable, as proven independently in \cite{Athanassenas, Vogel} for $n=3$, and in \cite{Pedrosa-Ritore} for $3\leq n\leq 8.$ In fact, when $n\leq 8$ the infinite array of balls and the cylinders are the only stable periodic solutions to the constant mean curvature equation; see \cite[Corollary 5.4]{Pedrosa-Ritore}.

When $n\geq 10$ it has been proven in \cite[Proposition 3.4]{Pedrosa-Ritore} that there exist unduloids which are minimizers under certain volume constraints (namely for volumes which are close to that of a ball that has the diameter of one period). For large volumes, however, cylinders are the only minimizers in all dimensions, as  shown in \cite[Theorem 1.1]{Ritore-Vernadikis}. For $n=9$, \cite[Section 5]{Koiso} provides numerical evidence that there exist stable unduloids in this dimension (as it is the case for $n\geq 10$).

A simpler fact is to characterize the stability of straight cylinders, still in the local case. A cylinder $\{(x_1,x')\in\R^n:\,|x'|<R\}$ is stable if and only if
\begin{equation}\label{stab-cyl} 
  R^2\geq n-2,
\end{equation}
as can be verified by a simple calculation (see \cite{Athanassenas} in the case $n=3$ and use the Euler-Lagrange equation \cite[Equation (4)]{Koiso} for general $n$). 

In view of these facts for the local case, and at least in low dimensions, one may expect the minimizers of \Cref{scan thm1} to be the near-balls found in \cite{CFW2} for small volumes and straight cylinders for large volumes. To sustain this statement for small volumes (which is still an open problem), recall first that the function $u$ defining the minimizer in \eqref{Del frac min E u} could be identically zero in a subinterval $(a,\pi)$ of $(0,\pi)$ ---allowing the minimizer, in this way, to be a periodic array of near-balls. Recall also that \cite[Theorem 3]{Davila} proved that, for small volumes, a minimizer must be close (in a measure theoretical sense) to a periodic array of round balls. On the other hand, the periodic array of near-spheres found in \cite{CFW2} are obtained as small perturbations of a single round sphere, and hence remains as a stable configuration.

With regard to large volumes,  in \cite{CCM Stability} we will establish the following result on the stability of straight cylinders. Its conclusion supports the conjecture that the only minimizers for large volumes should be straight cylinders. By stability we mean that the second variation of $\PP_s[E_t]$ is nonnegative for all $2\pi$-periodic volume preserving variations $\{E_t\}_{t\in(-\epsilon,\epsilon)}$ of the set $E$.\footnote{The proof of Theorem~\ref{thm:intro:stability} in~\cite{CCM Stability} will rely mainly on the techniques developed in \cite{CFSW,CFW1}. In particular, it will involve rather lengthy technical computations. This is the reason why we present the proof of this third theorem in a separate paper.}

\begin{theorem}[\cite{CCM Stability}]
\label{thm:intro:stability}
Let $n\geq 2$ and $s\in (0,1)$.
 
Then, there exists a radius $R_s>0$ (depending only on $n$ and $s$) such that the set $E=\{(x_1,x')\in \R^n:\, |x'|<R\}$ is stable (for $\PP_s$ and among sets which are $2\pi$-periodic in the $x_1$ direction)  if and only if $R\geq R_s$. Moreover,
$$
  \lim_{s\uparrow 1}R^2_s=n-2.
$$
\end{theorem} 

Observe that the last statement fits with \eqref{stab-cyl}, also in the case $n=2$ ---a dimension for which all above questions are trivial in the local case,  but not in the nonlocal case. 
The proof of Theorem \ref{thm:intro:stability} in \cite{CCM Stability} uses tools from \cite{CFSW,CFW1} mainly, but also the symmetric decreasing periodic rearrangement of the present paper. Indeed, we will use the fact that, thanks to \Cref{thm:intro:symmetry of Delaunay}, to show stability we will only need to consider even variations of $E.$ In this way, we will be able to use the basis $\cos(k\cdot)$, $k\in\mathbb{Z},$ of even eigenfunctions to the linearized operator of $\PP_s$ at a straight cylinder $E.$

 \subsection{Plan of the paper}
\begin{itemize}

\item[---]
In Section \ref{s delaunay} we prove \Cref{scan thm1}, even though we will borrow the results of  \Cref{section:rearrangement} only to prove the symmetry and monotonicity properties of minimizers claimed in the equality \eqref{Del frac min E u} of $(i)$ and in $(ii)$.

\item[---] Section \ref{section:rearrangement} deals with the behaviour of $\PP_s$ under cylindrical rearrangement and under symmetric decreasing periodic rearrangement. In particular, here we prove Theorem~\ref{thm:intro:symmetry of Delaunay}.

\item[---] In \Cref{section:periodic min in R2} we present a more explicit expression for $\PP_s$ when
 $n=2$. This gives an alternative proof of the existence part in Theorem \ref{scan thm1} for $n=2$.

\item[---] In \Cref{Appendix:Heat Kernel} we provide a simple proof of a well-known monotonicity property of the fundamental solution of the one-dimensional heat equation under periodic boundary conditions. The result is crucially used in the proof of \Cref{thm:intro:symmetry of Delaunay}.

\end{itemize}


\section{Existence of periodic minimizers}\label{s delaunay}
In this section we establish Theorem \ref{scan thm1}. 
Before stating its key lemma, we first prove an elementary one which will be used several times in the paper. 

\begin{lemma}\label{gy:lemma:interch integ periodic}
Let $m$ and $l$ be two positive integers, $u=(u_1,\ldots,u_m)$ and $v=(v_1,\ldots,v_l)$ be $2\pi$-periodic functions in $\R$ with values in $\R^m$ and $\R^l$, respectively. Suppose that $\Phi:\R^m\times\R^l\times[0,+\infty)\to\R$ is such that the function $(x,y)\mapsto \Phi(u(x),v(y),|x-y|)$ belongs to $L^1((-\pi,\pi)\times\R).$

Then, this last function belongs also to $L^1(\R\times(-\pi,\pi))$ and
\begin{equation}\label{gy:eq:interch integ per}
  \int_{-\pi}^{\pi}\!dx\int_{\R}\!dy\;\Phi\big(u(x),v(y),|x-y|\big)=
  \int_{\R}\!dx\int_{-\pi}^{\pi}\!dy\;\Phi\big(u(x),v(y),|x-y|\big).
\end{equation}
\end{lemma}

\begin{proof}
In the integral on the left-hand side of \eqref{gy:eq:interch integ per}, both for $|\Phi|$ and for $\Phi$, we can write the domain of integration $\R$ as the union over all $\big((2k-1)\pi,(2k+1)\pi\big),$ where $k\in \Z$. Thus the left-hand side is the countable sum of integrals over the domains $(x,y)\in (-\pi,\pi)\times \big((2k-1)\pi,(2k+1)\pi\big)$.  For each $k\in \Z$ we carry out the change of variables $x=\overline{x}+2k\pi$ and $y=\overline{y}+2k\pi, $ where $(\overline{x},\overline{y})\in \big((-2k-1)\pi,(-2k+1)\pi\big)\times (-\pi,\pi). $ Now, the periodicity of $u$ and $v$, and the identity $|x-y|=|\overline{x}-\overline{y}|$, yield the claim of the lemma.
\end{proof}

The following lemma will be a key ingredient in the proof of existence of minimizer.

\begin{lemma}\label{scan lemma 2}
Let $n\geq 2$ and $E\subset\R^n$ be a set of the form
\begin{equation}
E=\{x\in\R^n:\,|x'|<u(x_1)\}
\end{equation}
with $u:\R\to\R$ measurable, nonnegative, and $2\pi$-periodic in $x_1$. 

Then,
\begin{equation}\label{p2}
\PP_s(u):=\PP_s[E] = \int_{-\pi}^{\pi}\!dx_1\int_{-\infty}^{+\infty}\!dy_1\,
{|x_1-y_1|^{n-2-s}}\,\phi\Big(\frac{u(x_1)}{|x_1-y_1|},\frac{u(y_1)}{|x_1-y_1|}\Big),
\end{equation}
where, for nonnegative numbers $p$ and $q$, 
\begin{equation}\label{p3}
\phi(p,q):=\int_{\{w'\in\R^{n-1}:\,|w'|<p\}}\!dw'\int_{\{z'\in\R^{n-1}:\,|z'|>q\}}\!dz'
\,(1+|w'-z'|^2)^{-\frac{n+s}{2}}.
\end{equation}

Moreover, we have
\begin{equation}\label{p4}
\phi(p,q)\geq c_0(p^{n-1}-q^{n-1})-2c_0\qquad\text{if }0\leq q\leq p,
\end{equation}
for some constant $c_0>0$ depending only on $n$.
\end{lemma}

\begin{proof}
From \eqref{p}, we have
\begin{equation}
\begin{split}
\PP_s[E] &= \int_{-\pi}^{\pi}\!dx_1\int_{-\infty}^{+\infty}\!dy_1
\int_{\{x'\in\R^{n-1}:\,|x'|<u(x_1)\}}\!dx'\\
&\hskip98pt\int_{\{y'\in\R^{n-1}:\,|y'|>u(y_1)\}}\!dy'\,
(|x_1-y_1|^2+|x'-y'|^2)^{-\frac{n+s}{2}}\\
&= \int_{-\pi}^{\pi}\!dx_1\int_{-\infty}^{+\infty}\!dy_1
\,|x_1-y_1|^{-n-s}
\int_{\{x'\in\R^{n-1}:\,|x'|<u(x_1)\}}\!dx'\\
&\hskip98pt\int_{\{y'\in\R^{n-1}:\,|y'|>u(y_1)\}}\!dy'\,
\Big(1+\Big|\frac{x'}{|x_1-y_1|}-\frac{y'}{|x_1-y_1|}\Big|^2\Big)^{-\frac{n+s}{2}}\\
&= \int_{-\pi}^{\pi}\!dx_1\int_{-\infty}^{+\infty}\!dy_1
\,|x_1-y_1|^{n-2-s}
\int_{\big\{w'\in\R^{n-1}:\,|w'|<\frac{u(x_1)}{|x_1-y_1|}\big\}}\!dw'\\
&\hskip98pt\int_{\big\{z'\in\R^{n-1}:\,|z'|>\frac{u(y_1)}{|x_1-y_1|}\big\}}\!dz'\,
(1+|w'-z'|^2)^{-\frac{n+s}{2}},
\end{split}
\end{equation}
which gives \eqref{p2}-\eqref{p3}.

Next, to prove \eqref{p4} for $0\leq q\leq p$, we see that
\begin{equation}\label{p5}
\begin{split}
\phi(p,q)&=\int_{\{w'\in\R^{n-1}:\,|w'|<p\}}\!dw'
\int_q^{+\infty}\!dr\,r^{n-2}\int_{\Sph^{n-2}}\!d\mathcal{H}^{n-2}(\sigma)\,
(1+|w'-r\sigma|^2)^{-\frac{n+s}{2}}\\
&=\int_q^{+\infty}\!dr\,r^{n-2}\int_{\Sph^{n-2}}\!d\mathcal{H}^{n-2}(\sigma)\,
\int_{\{w'\in\R^{n-1}:\,|w'|<p\}}\!dw'\,(1+|w'-r\sigma|^2)^{-\frac{n+s}{2}}\\
&=|\Sph^{n-2}|\int_q^{+\infty}\!dr\,r^{n-2}
\int_{\{w'\in\R^{n-1}:\,|w'|<p\}}\!dw'\,(1+|w'-re_1'|^2)^{-\frac{n+s}{2}}\\
&\geq|\Sph^{n-2}|\int_q^{p}\!dr\,r^{n-2}
\int_{\{w'\in\R^{n-1}:\,|w'|<p\}}\!dw'\,(1+|w'-re_1'|^2)^{-\frac{n+s}{2}},
\end{split}
\end{equation}
where $e_1'=(1,0,\ldots,0)\in\R^{n-1}$.

We now claim that, for some constant $c>0$ depending only on $n$,\footnote{Here and throughout the paper, $B_r^{k}(y)$ denotes the open ball in $\R^k$ of radius $r$ centered at $y\in\R^k$. We also use the notation $B_r^{k}$ when $y=0$, and $B_r$ if in addition $k=n$.}
\begin{equation}\label{p6}
|B^{n-1}_1(re_1')\cap B^{n-1}_p|\geq c\qquad\text{if }1\leq r\leq p.
\end{equation}
This is simple. We have
$B^{n-1}_p\supset B^{n-1}_r\supset B^{n-1}_1((r-1)e_1').$ Hence,
\begin{equation}
  |B^{n-1}_1(re_1')\cap B^{n-1}_p|\geq |B^{n-1}_1(re_1')\cap B^{n-1}_1((r-1)e_1')|
  =|B^{n-1}_1\cap B^{n-1}_1(-e_1')|=:c.
\end{equation}

Therefore, from \eqref{p5} and \eqref{p6} we deduce that, for $1\leq q\leq p$,
\begin{equation}
\begin{split}
\phi(p,q)&\geq|\Sph^{n-2}|\int_q^{p}\!dr\,r^{n-2}
\int_{B^{n-1}_1(re_1')\cap B^{n-1}_p}\!dw'\,(1+|w'-re_1'|^2)^{-\frac{n+s}{2}}\\
&\geq2^{-\frac{n+s}{2}}|\Sph^{n-2}|\int_q^{p}\!dr\,r^{n-2}|B^{n-1}_1(re_1')\cap B^{n-1}_p|\\
&\geq2^{-\frac{n+1}{2}}|\Sph^{n-2}|c\int_q^{p}\!dr\,r^{n-2}
=:2c_0(p^{n-1}-q^{n-1}).
\end{split}
\end{equation}
We conclude that \eqref{p4} holds if $1\leq q\leq p$.

Next, if $q<1$ and $p^{n-1}\geq 2$, we see that $\phi(p,q)\geq\phi(p,1)
\geq 2c_0(p^{n-1}-1)$ by definition of $\phi$ and the previous bound. But $p^{n-1}-1\geq(p^{n-1}-q^{n-1})/2$ since we assumed $p^{n-1}\geq 2$. Thus, in this case \eqref{p4} also holds.

It only remains to consider the case $0\leq q<1$ and $p^{n-1}< 2$. But then we have
\begin{equation}
c_0(p^{n-1}-q^{n-1})-2c_0\leq c_0 p^{n-1}-2c_0<0\leq\phi(p,q),
\end{equation}
which finishes the proof.
\end{proof}

We now establish Theorem \ref{scan thm1}. 
The proof of existence of a minimizer relies only on the two previous basic lemmas and on the fact (already used in \cite{Davila}) that $\PP_s$ does not increase under cylindrical rearrangement (as stated in \Cref{gy:lemma:cylindrical rearr}). Instead, the proof of the eveness and nonincreasing property of any minimizer uses our new result \Cref{thm:intro:symmetry of Delaunay} on the periodic rearrangement.

\begin{proof}[Proof of Theorem~\ref{scan thm1}]
We split the proof into four parts.

\smallskip

{\it \underline{Part 1} $($\underline{Existence of minimizer}$)$.}
Here we prove the existence of a minimizer. Note first that $\PP_s\geq 0$ and that $\PP_s\not \equiv +\infty$ ---this last fact will be easily seen in Part 4 of the proof when we check that $\PP_s$ is finite when computed on a periodic array of balls. 

Given $\mu>0,$ we take a minimizing sequence of sets with volume $\mu$ within the slab $\{-\pi<x_1<\pi\}$. In view of Lemma~\ref{gy:lemma:cylindrical rearr}, the sets may be assumed to be of the form
\begin{equation}
E_k=\{x\in\R^n:\, |x'|<u_k(x_1)\}
\end{equation}
for some measurable, nonnegative, and $2\pi$-periodic functions $u_k$. 

By Lemma~\ref{scan lemma 2} we have that, for the constant $c_0$ in \eqref{p4} and some positive constant $C_{n,s}$ depending only on $n$ and $s$,
\begin{equation}
\begin{split}
\int_{-\pi}^{\pi}\!dx_1&\int_{-\pi}^{\pi}\!dy_1\,{|x_1-y_1|^{n-2-s}}
\Big\{c_0+\phi\Big(\frac{u_k(x_1)}{|x_1-y_1|},\frac{u_k(y_1)}{|x_1-y_1|}\Big)\Big\}\\
&= C_{n,s}\,c_0+\int_{-\pi}^{\pi}\!dx_1\int_{-\pi}^{\pi}\!dy_1\,{|x_1-y_1|^{n-2-s}}
\phi\Big(\frac{u_k(x_1)}{|x_1-y_1|},\frac{u_k(y_1)}{|x_1-y_1|}\Big)\\
&\leq C_{n,s}\,c_0+\PP_s[E_k]\leq C
\end{split}
\end{equation}
since $\phi\geq 0$, where $C$ is a constant independent of $k$.
From this, symmetrizing in $x_1$ and $y_1$, using \eqref{p4} and, again, $\phi\geq 0$, we see that 
\begin{equation}\label{gy:v12 star eq} 
\begin{split}
C&\geq \int_{-\pi}^{\pi}\!dx_1\int_{-\pi}^{\pi}\!dy_1\,{|x_1-y_1|^{n-2-s}}
\Big\{c_0+\phi\Big(\frac{u_k(x_1)}{|x_1-y_1|},\frac{u_k(y_1)}{|x_1-y_1|}\Big)\Big\}\\
&=\frac{1}{2}\int_{-\pi}^{\pi}\!dx_1\int_{-\pi}^{\pi}\!dy_1\,{|x_1-y_1|^{n-2-s}}
\Big\{c_0+\phi\Big(\frac{u_k(x_1)}{|x_1-y_1|},\frac{u_k(y_1)}{|x_1-y_1|}\Big)\Big\}\\
&\hskip80pt+\frac{1}{2}\int_{-\pi}^{\pi}\!dx_1\int_{-\pi}^{\pi}\!dy_1\,{|x_1-y_1|^{n-2-s}}
\Big\{c_0+\phi\Big(\frac{u_k(y_1)}{|x_1-y_1|},\frac{u_k(x_1)}{|x_1-y_1|}\Big)\Big\}\\
&\geq\frac{1}{2}\int_{-\pi}^{\pi}\!dx_1\int_{-\pi}^{\pi}\!dy_1\,{|x_1-y_1|^{n-2-s}}\\
&\hskip80pt\times\Big\{2c_0+\phi\Big(\frac{\max\{u_k(x_1),u_k(y_1)\}}{|x_1-y_1|},
\frac{\min\{u_k(x_1),u_k(y_1)\}}{|x_1-y_1|}\Big)\Big\}\\
&\geq\frac{c_0}{2}\int_{-\pi}^{\pi}\!dx_1\int_{-\pi}^{\pi}\!dy_1\,{|x_1-y_1|^{n-2-s}}
\frac{|u_k(x_1)^{n-1}-u_k(y_1)^{n-1}|}{|x_1-y_1|^{n-1}}\\
&=\frac{c_0}{2}\int_{-\pi}^{\pi}\!dx_1\int_{-\pi}^{\pi}\!dy_1\, 
\frac{|u_k(x_1)^{n-1}-u_k(y_1)^{n-1}|}{|x_1-y_1|^{1+s}}.
\end{split}
\end{equation}

This gives a uniform bound on the 
$W^{s,1}(-\pi,\pi)$ seminorm of the functions $u_k^{n-1}$. At the same time, thanks to our volume constraint, we know that
$|B_1^{n-1}|\int_{-\pi}^\pi dx_1\, u_k(x_1)^{n-1}=\mu$ for all~$k$.
Therefore, by the fractional Sobolev compactness theorem (see \cite[Theorem 7.1]{Valdinoci1}), a subsequence of $\{u_k^{n-1}\}$ converges strongly in $L^1(-\pi,\pi)$ to a nonnegative function
$v\in W^{s,1}(-\pi,\pi)$. As a consequence, if we set $u:=v^{1/(n-1)}$, then 
$|B_1^{n-1}|\int_{-\pi}^\pi dx_1\, u(x_1)^{n-1}=\mu$.
Thus, extending $u$ to be $2\pi$-periodic in $\R$, $u$ is an admissible competitor. Finally, since the function $\phi$ in \eqref{p2} is nonnegative, Fatou's lemma gives that $E:=\{x\in\R^n:\, |x'|<u(x_1)\}$ is a minimizer. 

The fact that $v=u^{n-1}\in W^{s,1}(-\pi,\pi)$ yields $u\in W^{s/(n-1),n-1}(-\pi,\pi)$, by simply using that $|t-s|^{n-1}\le |t^{n-1} - s^{n-1}|$ for all nonnegative numbers $t$ and $s$.

\smallskip

{\it\underline{Part 2} $($\underline{Symmetry of minimizers}$).$}
To prove $(i)$ in \Cref{scan thm1}, let $E$ be any minimizer. As in Part 1 above, by \Cref{gy:lemma:cylindrical rearr} we must have $E=\{|x'|<u(x_1)\}$ for some nonnegative $2\pi$-periodic function $u$. By \eqref{gy:v12 star eq} (with $u_k$ replaced by $u$), we see that $u^{n-1}\in W^{s,1}(-\pi,\pi)$. As before, this leads to $u\in W^{s/(n-1),n-1}(-\pi,\pi)$.

Next, statement $(ii)$ in the theorem follows immediately from \Cref{thm:intro:symmetry of Delaunay}.

\smallskip

{\it\underline{Part 3}  $($\underline{Euler-Lagrange equation}$)$.} In this part we address the proof of \Cref{scan thm1} $(iii)$. Firstly,
we will establish that if $E=\{x\in\R^n:\, |x'|<u(x_1)\}$ is a minimizer, $x_1^0\in\R$ is such that $u$ is $C^{1,\alpha}$ in a neighborhood $U(x_1^0)$ of $x_1^0$ for some $\alpha>s$, and $u(x_1^0)>0$, then the nonlocal mean curvature of $\partial E$ is constant in a neighborhood of $x_1^0$. Afterwards, we will show that the constant does not depend on the point $x_1^0$ in whose neighborhood $u$ is $C^{1,\alpha}$.  

Note first that we may assume $x_1^0\in(-\pi,\pi)$, by the periodicity of $E$.
We take any smooth function $\xi=\xi(x_1)$ with compact support in the neighborhood $U(x_1^0)$ of $x_1^0$, small enough to be contained in $(-\pi,\pi)$, and we extend $\xi$ to be $2\pi$-periodic. The first variation of the volume functional ---that is, the derivative of 
$|B_1^{n-1}|\int_{-\pi}^\pi dx_1\,  (u(x_1)+\epsilon \xi(x_1))^{n-1}$ with respect to $\epsilon$ evaluated at $\epsilon=0$---
is
\begin{equation}\label{p10}
(n-1)|B_1^{n-1}|\int_{-\pi}^\pi dx_1\,  u(x_1)^{n-2}\xi(x_1),
\end{equation}
while that of $\PP_s$ is 
\begin{equation}\label{p11}
\begin{split}
\frac{d}{d\epsilon}\biggr\rvert_{\epsilon=0}\PP_s(u+\epsilon\xi)
=\int_{-\pi}^{\pi}\!dx_1 
  \int_{-\infty}^{+\infty}\!dy_1
  &
  \,{|x_1-y_1|^{n-2-s}}
\Big\{\phi_p\Big(\frac{u(x_1)}{|x_1-y_1|},\frac{u(y_1)}{|x_1-y_1|}\Big)
\frac{\xi(x_1)}{|x_1-y_1|}\\
&\phantom{aaaaaaaaaa}
+\phi_q\Big(\frac{u(x_1)}{|x_1-y_1|},\frac{u(y_1)}{|x_1-y_1|}\Big)
\frac{\xi(y_1)}{|x_1-y_1|}\Big\}.
\end{split}
\end{equation}
We next see that this expression leads to the nonlocal mean curvature of $\partial E$ at $x=(x_1,x')$. In what follows, all expressions can be checked to be well defined (as for the nonlocal mean curvature of general, nonperiodic, $C^{1,\alpha}$ sets) using that $u$ is $C^{1,\alpha}$, with $\alpha>s$, in the neighborhood of $x_1^0$ where $\xi$ has its support.

Since, by \Cref{gy:lemma:interch integ periodic} (used with $m=1$ and $l=2$),
\begin{equation}
\begin{split}
  \int_{-\pi}^{\pi}\!dx_1&\int_{-\infty}^{+\infty}\!dy_1\,\frac{\xi(y_1)}{|x_1-   
  y_1|^{3+s-n}}\,
   \phi_q\Big(\frac{u(x_1)}{|x_1-y_1|},\frac{u(y_1)}{|x_1-y_1|}\Big)
   \\
   &=
  \int_{-\infty}^{+\infty}\!dx_1\int_{-\pi}^{\pi}\!dy_1\,\frac{\xi(y_1)}{|x_1-   
  y_1|^{3+s-n}}\,
   \phi_q\Big(\frac{u(x_1)}{|x_1-y_1|},\frac{u(y_1)}{|x_1-y_1|}\Big)
\\
  &=\int_{-\pi}^{\pi}\!dx_1\int_{-\infty}^{+\infty}\!dy_1\,\frac{\xi(x_1)}{|x_1-y_1|^{3+s-n}}
\phi_q\Big(\frac{u(y_1)}{|x_1-y_1|},\frac{u(x_1)}{|x_1-y_1|}\Big),
\end{split}
\end{equation}
we deduce that \eqref{p11} is equal to
\begin{equation}
\begin{split}
\int_{-\pi}^{\pi}\!dx_1\int_{-\infty}^{+\infty}\!dy_1\,\frac{\xi(x_1)}{|x_1-y_1|^{3+s-n}}
\Big\{\phi_p\Big(\frac{u(x_1)}{|x_1-y_1|},\frac{u(y_1)}{|x_1-y_1|}\Big)
+\phi_q\Big(\frac{u(y_1)}{|x_1-y_1|},\frac{u(x_1)}{|x_1-y_1|}\Big)
\Big\}.
\end{split}
\end{equation}
Thus, by definition of $\phi$, this last expression becomes 
\begin{equation}
\begin{split}
\int_{-\pi}^{\pi}\!d&x_1\,\xi(x_1)\int_{-\infty}^{+\infty}\!dy_1\,{|x_1-y_1|^{n-3-s}}
\Big\{\phi_p\Big(\frac{u(x_1)}{|x_1-y_1|},\frac{u(y_1)}{|x_1-y_1|}\Big)
+\phi_q\Big(\frac{u(y_1)}{|x_1-y_1|},\frac{u(x_1)}{|x_1-y_1|}\Big)
\Big\}\\
&=\int_{-\pi}^{\pi}\!dx_1\,\xi(x_1)\int_{-\infty}^{+\infty}\!dy_1\,{|x_1-y_1|^{n-3-s}}\\
&\hskip25pt\times
\Big\{\int_{\big\{w'\in\R^{n-1}:\,|w'|=\frac{u(x_1)}{|x_1-y_1|}\big\}}\!dw'
\int_{\big\{z'\in\R^{n-1}:\,|z'|>\frac{u(y_1)}{|x_1-y_1|}\big\}}\!dz'\,
(1+|w'-z'|^2)^{-\frac{n+s}{2}}\\
&\hskip35pt-\int_{\big\{w'\in\R^{n-1}:\,|w'|<\frac{u(y_1)}{|x_1-y_1|}\big\}}\!dw'
\int_{\big\{z'\in\R^{n-1}:\,|z'|=\frac{u(x_1)}{|x_1-y_1|}\big\}}\!dz'\,
(1+|w'-z'|^2)^{-\frac{n+s}{2}}\Big\}\\
&=|\Sph^{n-2}|\int_{-\pi}^{\pi}\!dx_1\xi(x_1)u(x_1)^{n-2}\int_{-\infty}^{+\infty}\!dy_1\,{|x_1-y_1|^{n-3-s}}
\Big(\frac{1}{|x_1-y_1|}\Big)^{n-2}\\
&\hskip57pt\times
\Big\{\int_{\big\{z'\in\R^{n-1}:\,|z'|>\frac{u(y_1)}{|x_1-y_1|}\big\}}\!dz'\,
\Big(1+\Big|z'-\frac{u(x_1)e_1'}{|x_1-y_1|}\Big|^2\Big)^{-\frac{n+s}{2}}\\
&\hskip67pt-\int_{\big\{w'\in\R^{n-1}:\,|w'|<\frac{u(y_1)}{|x_1-y_1|}\big\}}\!dw'
\Big(1+\Big|w'-\frac{u(x_1)e_1'}{|x_1-y_1|}\Big|^2\Big)^{-\frac{n+s}{2}}\Big\}.
\end{split}
\end{equation}
The integral on $dy_1$ over $(-\infty,+\infty)$ is equal to 
\begin{equation}
\begin{split}
\int_{-\infty}^{+\infty}&\!dy_1\,{|x_1-y_1|^{-1-s}}\\
&\hskip10pt\times
\Big\{\int_{\{\overline z'\in\R^{n-1}:\,|\overline z'|>{u(y_1)}\}}\!d\overline z'\,
|x_1-y_1|^{1-n}\Big(1+\Big|\frac{\overline z'}{|x_1-y_1|}-\frac{u(x_1)e_1'}{|x_1-y_1|} \Big|^2\Big)^{-\frac{n+s}{2}}\\
&\hskip30pt-\int_{\{\overline w'\in\R^{n-1}:\,|\overline w'|<{u(y_1)}\}}\!d\overline w'\,
|x_1-y_1|^{1-n}\Big(1+\Big|\frac{\overline w'}{|x_1-y_1|}-\frac{u(x_1)e_1'}{|x_1-y_1|} \Big|^2\Big)^{-\frac{n+s}{2}}\Big\}\\
=&\int_{-\infty}^{+\infty}\!dy_1\,\Big\{
\int_{\{\overline z'\in\R^{n-1}:\,|\overline z'|>{u(y_1)}\}}\!d\overline z'\,
(|x_1-y_1|^2+|u(x_1)e_1'-\overline z'|^2)^{-\frac{n+s}{2}}\\
&\hskip55pt-\int_{\{\overline w'\in\R^{n-1}:\,|\overline w'|<{u(y_1)}\}}\!d\overline w'\,
(|x_1-y_1|^2+|u(x_1)e_1'-\overline w'|^2)^{-\frac{n+s}{2}}
\Big\}.
\end{split}
\end{equation}

This last expression is, by definition, the nonlocal mean curvature of $\partial E$ at the point $x=(x_1,u(x_1)e_1')=(x_1,u(x_1),0,\ldots,0)$, which is defined by
\begin{equation}
H_{s}[E](x) = \int_{\R^n}dy\,\frac{\chi_{E^c}(y)-\chi_{E}(y)}{|x-y|^{n+s}}.
\end{equation}
Thus, by \eqref{p11}, we have shown that
\begin{equation}
\label{gy:eq:1st var of Ps}
  \frac{d}{d\epsilon}\biggr\rvert_{\epsilon=0}\PP_s(u+\epsilon\xi)
=|\Sph^{n-2}|\int_{-\pi}^{\pi}\!dx_1\, \xi(x_1) u(x_1)^{n-2}
H_s[E](x_1,u(x_1)e_1').
\end{equation}
From this, \eqref{p10}, and the Lagrange multiplier 
rule,\footnote{Note that we are actually using the Lagrange multiplier rule not for the original functional $\PP_s$ but for another functional $\overline{\PP}_s$ given by the restriction of $\PP_s$ to functions  which agree with $u$ outside of the set where $u$ is regular. The Lagrange multiplier rule applies to $\overline{\PP}_s$, since $u$ is also a minimizer for $\overline{\PP}_s$ and $\overline{\PP}_s$ is differentiable.}
we see that, for some constant $\lambda\in\R$,
$$
  H_s[E](x_1,u(x_1)e_1')=\lambda \quad\text{ for all }x_1\in U(x_1^0).
$$

It remains to show that $\lambda$ does not depend on $x_1^0$. For this purpose, let $U(x_1^0)$ and $U(y_1^0)$ be two neighborhoods of two different points $x^0_1$ and $y^0_1$ where $u$ is regular. Now, we can simply repeat the previous argument by taking this time a smooth function $\xi$ with compact support in the union $U(x_1^0)\cup U(y_1^0)$. Taking such variations $\xi$ leads to the conclusion that $H_s[E](x_1,u(x_1)e_1')$ is constant in $U(x_1^0)\cup U(y_1^0)$. This proves \Cref{scan thm1} $(iii)$.

\vspace{0.2cm}

{\it\underline{Part 4} $($\underline{Small volumes}$).$} We finally address the proof of \Cref{scan thm1} $(iv)$. That is, we prove (uniformly as $s\uparrow1$ and $s\downarrow0$) that if the volume constraint $\mu$ is small enough, any minimizer $E$ of $\PP_s$ is not a straight cylinder. 
To do so, we simply use an energy comparison with the configuration $\BB_r$ given by a periodic array of disjoint balls
\begin{equation}\label{p7}
\BB_r:=\{x\in\R^n:\,|x-2k\pi e_1|<r\text{ for some }k\in\Z\},
\end{equation}
where $e_1=(1,0,\ldots,0)\in\R^{n}$ and $r<\pi$. 
To satisfy the volume constraint, the radius $r=r(\mu)$ of the balls must be taken to satisfy $|B_1| \, r(\mu)^n=\mu$.

Now, by expression \eqref{p} we have, for different constants $C_n$ depending only on $n$, 
\begin{equation}\label{p12}
\begin{split}
\PP_s[\BB_r] &= \int_{B_r}\!dx\int_{\R^n\setminus \BB_r} \frac{dy}{|x-y|^{n+s}}\\
&\leq \int_{B_r}\!dx\int_{\R^n\setminus B_r} \frac{dy}{|x-y|^{n+s}}
= r^{n-s}\int_{B_1}\!dx\int_{\R^n\setminus B_1} \frac{dy}{|x-y|^{n+s}}\\
&\leq r^{n-s}\int_{B_1}\!dx\int_{\R^n\setminus B_{1-|x|}(x)} \frac{dy}{|x-y|^{n+s}}
\leq C_n\, r^{n-s}\int_{B_1}\!dx\int_{1-|x|}^{+\infty}\!dt\,t^{-1-s}\\
&=\frac{C_n}{s}\,r^{n-s}\int_{B_1}\!dx\,(1-|x|)^{-s}
=\frac{C_n}{s}\,r^{n-s}|\Sph^{n-1}|\int_0^1\!d\rho\,\rho^{n-1}(1-\rho)^{-s}
\\
&\leq\frac{C_n}{s(1-s)}\,r^{n-s}=  \frac{C_n}{s(1-s)}\,\mu^{\frac{n-s}{n}}.
\end{split}
\end{equation}

Instead, for the straight cylinder satisfying the volume constraint, which is given by 
$E_\mu:=\{x\in\R^n:\,|x'|<c_1\,\mu^{1/(n-1)}\}$ for some constant $c_1$ depending only on $n$, we have
\begin{equation}
\begin{split}
\PP_s[E_\mu] &= \int_{\{x\in\R^{n}:\,|x_1|<\pi,\,|x'|<c_1\,\mu^{1/(n-1)}\}}\!dx\,
\int_{\{y\in\R^n:\,|y'|>c_1\,\mu^{1/(n-1)}\}} \frac{dy}{|x-y|^{n+s}}\\
&=c_1^{n-s}\,\mu^{\frac{n-s}{n-1}}
\int_{\big\{\overline x\in\R^{n}:\,|\overline x_1|<\frac{\pi}{c_1\,\mu^{1/(n-1)}}
,\,|\overline x'|<1\big\}}\!d\overline x\,
\int_{\{\overline y\in\R^n:\,|\overline y'|>1\}} \frac{d\overline y}{|\overline x-\overline y|^{n+s}}\\
&=c_1^{n-s}\,\mu^{\frac{n-s}{n-1}}\frac{2\pi}{c_1\,\mu^{1/(n-1)}}
\int_{-\infty}^{+\infty}\!dt\int_{\{\overline x'\in\R^{n-1}:\,|\overline x'|<1\}}\!d\overline x'\!
\int_{\{\overline y'\in\R^{n-1}:\,|\overline y'|>1\}}\!d\overline y'\, 
{(t^2+|\overline x'-\overline y'|^2)^{-\frac{n+s}{2}}}\\
&\geq c_1^{n-s-1}\,\mu^{\frac{n-1-s}{n-1}}2\pi \, 2^{-\frac{n+s}{2}}
\int_{\{\overline x'\in\R^{n-1}:\,|\overline x'|<1\}}\!d\overline x'
\int_{\{\overline y'\in\R^{n-1}:\,|\overline y'|>1\}}\!d\overline y'
\int_{-|\overline x'-\overline y'|}^{|\overline x'-\overline y'|}\!dt\, 
|\overline x'-\overline y'|^{-n-s}\\
&\geq c_1^{n-s-1}\,\mu^{\frac{n-1-s}{n-1}}2\pi\,2^{-\frac{n+s}{2}}2
\int_{\{\overline x'\in\R^{n-1}:\,|\overline x'|<1\}}\!d\overline x'
\int_{\{\overline y'\in\R^{n-1}:\,|\overline y'|>1\}}\!d\overline y'\, 
|\overline x'-\overline y'|^{1-n-s}.
\end{split}
\end{equation}
To bound the last double integral from below, we denote
$\overline y'=z'=(z'_1,z'')\in\R\times\R^{n-2}$ to get
\begin{equation}
\begin{split}
\int_{B_1^{n-1}}\!d\overline x'&
\int_{\R^{n-1}\setminus B_1^{n-1}}\!d\overline y'\, 
|\overline x'-\overline y'|^{1-n-s}
\\
&=\big|\Sph^{n-2}\big|\int_{0}^1\!d \rho\,\rho^{n-2}
\int_{\R^{n-1}\setminus B_1^{n-1}}\!dz'\, 
|\rho\, e_1'-z'|^{1-n-s}
\\
&\geq \big|\Sph^{n-2}\big|\int_{0}^1\!d \rho\,\rho^{n-2}
\int_{\{ z\in\R^{n-1}:\,z'_1>1\}}\!dz'\, 
|\rho\, e_1'-z'|^{1-n-s}
\\
&=\big|\Sph^{n-2}\big|\int_{0}^1\!d \rho\,\rho^{n-2}
\int_1^{+\infty}\!dz'_1\int_{\R^{n-2}}\!dz''
\,\big( (\rho-z'_1)^2+|z''|^2\big)^{\frac{1-n-s}{2}}
\\
&\geq \big|\Sph^{n-2}\big|\int_{0}^1\!d \rho\,\rho^{n-2}
\int_1^{+\infty}\!\!dz'_1\!\int_{\{z''\in\R^{n-2}:\,|z''|<|\rho-z_1'|\}}
\!\!\!\!\!\!\!\!\!\!dz''\,\big( (\rho-z'_1)^2+|z''|^2\big)^{\frac{1-n-s}{2}}
\\
&\geq \big|\Sph^{n-2}\big|
2^{\frac{1-n-s}{2}}\int_{0}^1\!d \rho\,\rho^{n-2}
\int_1^{+\infty}\!dz'_1\int_{\{z''\in\R^{n-2}:\,|z''|<|\rho-z'_1|\}}
\!\!\!\!\!\!\!\!\!\!dz''\,|\rho-z'_1|^{1-n-s}\\
&=\big|\Sph^{n-2}\big|2^{\frac{1-n-s}{2}} \big|\Sph^{n-3}\big|
\int_{0}^1\!d \rho\,\rho^{n-2}
\int_1^{+\infty}\!dz'_1\,|\rho-z'_1|^{-1-s}\\
&=\frac{1}{s} \big|\Sph^{n-2}\big|2^{\frac{1-n-s}{2}} \big|\Sph^{n-3}\big|
\int_{0}^1d \rho\,\rho^{n-2}(1-\rho)^{-s}.
\end{split}
\end{equation}
Here we implicitly assumed that $n\geq 3$. However, when $n=2$ we do not need to introduce the integral over $z''\in\R^{n-2}$, and the same bound (with the constant $\frac{1}{s} \big|\Sph^{n-2}\big|2^{\frac{1-n-s}{2}} \big|\Sph^{n-3}\big|$ replaced by $\frac{2}{s}$) follows in an easier way. 
Observe now that
\begin{equation}
\begin{split}
\int_{0}^1\!\rho^{n-2}(1-\rho)^{-s}\,d \rho
\geq 2^{2-n}\int_{1/2}^1(1-\rho)^{-s}\, d\rho
=\frac{2^{1-n+s}}{1-s}\geq\frac{2^{1-n}}{1-s}.
\end{split}
\end{equation}
Hence, since $2^{-s}\geq 2^{-1}$ and $c_1^{-s}\geq \min(1,c_1^{-1})$, we deduce that
\begin{equation}\label{p13 aux 2}
\begin{split}
\PP_s[E_\mu] 
\geq\frac{c_n}{s(1-s)}\mu^{\frac{n-1-s}{n-1}}
\end{split}
\end{equation}
for some constant $c_n>0$ depending only on $n$. 

From this and \eqref{p12}, we conclude that $\PP_s[E_\mu]>\PP_s[\BB_{r(\mu)}]$ if $\mu$ is small enough, since clearly $(n-1-s)/(n-1)<(n-s)/n$.
\end{proof}


\section{Cylindrical and periodic symmetric decreasing rearrangements}
\label{section:rearrangement}


In this section we analyze the behavior of $\PP_s$ under the two rearrangements defined in Subsection~\ref{symmetry}, namely, the cylindrical rearrangement and the symmetric decreasing periodic rearrangement. 

In the previous section we have shown the existence of a cylindrically symmetric minimizer using the Riesz rearrangement inequality (the key ingredient in the proof of \Cref{gy:lemma:cylindrical rearr}).\footnote{The Riesz rearrangement inequality was already used  in \cite[Proposition 13]{Davila} for the functional $\PP_s$ acting on periodic sets, to conclude that the cylindrical rearrangement does not increase the functional $\PP_s$.}
In the next lemma we will see that characterizing the case of equality in the Riesz rearrangement inequality will actually yield that every minimizer is cylindrically symmetric. The same argument has been used in \cite[Lemma 4.2]{Malchiodi Novaga Pagliardini} to show the cylindrical symmetry of minimizers of a different functional.

Before stating the lemma we recall the notation introduced in Subsection~\ref{symmetry}. For a measurable set $E\subset\R^n$, given $x_1\in\R$ consider
\begin{equation}\label{def.rear.En1}
  E_{x_1}:=\{x'\in \R^{n-1}:\,(x_1,x')\in E\}.
\end{equation}
Then, denoting balls in $\R^{n-1}$ by $B_r^{n-1}$ and the $(n-1)$-dimensional Lebesgue measure
by $\mathcal{L}^{n-1}$, we define
\begin{equation}
E^{*\!\operatorname{cyl}}:=\left\{(x_1,x')\in \R^n:\,
  x'\in B_{r(x_1)}^{n-1},\,\text{ where }\mathcal{L}^{n-1}
  \big(B_{r(x_1)}^{n-1}\big)=\mathcal{L}^{n-1}(E_{x_1}),
  \right\},
\end{equation}
with the understanding that $B_{r(x_1)}^{n-1}=\R^{n-1}$ if $\mathcal{L}^{n-1}(E_{x_1})=+\infty$.
Notice that by definition
\begin{equation}
 \label{gy:eq:x1section and rearr interch}
   \left(E^{*\!\operatorname{cyl}}\right)_{x_1}=\left(E_{x_1}\right)^{*(n-1)},
\end{equation}
where $A^{*(n-1)}$ denotes the Schwarz rearrangement of a set $A\subset\R^{n-1},$ i.e., the ball in $\R^{n-1}$ centered at $0$ and of the same $(n-1)$-dimensional Lebesgue measure as $A$ ---in case that the measure of $A$ is infinite, we set $A^{*(n-1)}:=\R^{n-1}$.

The following lemma concerns our periodic fractional perimeter, but its statement does not require the periodicity of the set $E$.
Observe that $\PP_s[E]$ makes sense also for sets which are not periodic. Its proof will use the classical Riesz rearrangement inequality, as well as its strict version needed to discuss the case of equality.

\begin{lemma}
\label{gy:lemma:cylindrical rearr}
Let $E\subset\R^n$ be a measurable set such that 
$|E\cap((-\pi,\pi)\times\R^{n-1})|<+\infty$. 

Then, 
$\PP_s[E^{*\!\operatorname{cyl}}]\leq \PP_s[E]$. Moreover, if in addition $|E\cap((-\pi,\pi)\times\R^{n-1})|>0$ and $\PP_s[E]<+\infty$, the equality $\PP_s[E^{*\!\operatorname{cyl}}]= \PP_s[E]$ holds if and only if $E=E^{*\!\operatorname{cyl}}+(0,c)$ for some $c\in\R^{n-1}$ up to a set of measure zero. 
\end{lemma}

\begin{proof}
For simplicity of notation, given $x_1,y_1\in\R$ with $x_1\neq y_1$, define $g_{x_1,y_1}:\R^{n-1}\to \R$ by 
$$
  g_{x_1,y_1}(z'):=\left(|x_1-y_1|^2+|z'|^2\right)^{-\frac{n+s}{2}}.
$$
Note that $g_{x_1,y_1}\in L^1(\R^{n-1})$ is a positive and radially symmetric decreasing function.

To address the proof of 
the inequality $\PP_s[E^{*\!\operatorname{cyl}}]\leq \PP_s[E]$, let us first make some preliminary observations. Let $x_1,y_1\in\R$ with $x_1\neq y_1$, and assume that $\mathcal{L}^{n-1}(E_{x_1})<+\infty$. On the one hand, if $\mathcal{L}^{n-1}(E_{y_1})<+\infty$,
the Riesz rearrangement inequality \cite[Theorem 3.7]{Lieb Loss} yields that
\begin{equation}
 \label{gxy.eq.riesz}
  \int_{E_{x_1}}dx'\int_{E_{y_1}}dy'\,g_{x_1,y_1}(x'-y')
  \leq 
  \int_{(E_{x_1})^{*(n-1)}}dx'\int_{(E_{y_1})^{*(n-1)}}dy'\,g_{x_1,y_1}(x'-y').
\end{equation}
In addition, the strict rearrangement inequality \cite[Corollary 2.19]{Baernstein book} or 
\cite[Theorem 3.9]{Lieb Loss}\footnote{\cite[Theorem 3.9]{Lieb Loss} should include two hypotheses which are needed for its validity ---and which are explicitly stated in most of the results of \cite[Chapter 3]{Lieb Loss}. Namely, that the functions $f$ and $h$ on its statement must vanish at infinity and that neither $f$ nor $h$ is identically zero. Indeed, there is always equality in the Riesz rearrangement inequality if one of these two functions is constant.  
}
shows that the equality in \eqref{gxy.eq.riesz} holds (when $E_{x_1}$ and $E_{y_1}$ have finite measure) if and only if either 
$\mathcal{L}^{n-1}(E_{x_1})\mathcal{L}^{n-1}(E_{y_1})=0$ or
\begin{equation}\label{eq:case of eq in Riesz cylin}
E_{x_1}=B_{r(x_1)}^{n-1}+c(x_1,y_1)\quad\text{and}\quad E_{y_1}=B_{r(y_1)}^{n-1}+c(x_1,y_1)
\end{equation}
for some $r(x_1)$ and $r(y_1)$ belonging to $(0,+\infty)$ and for some
$c(x_1,y_1)\in\R^{n-1}$. 

On the other hand, if $\mathcal{L}^{n-1}(E_{y_1})=+\infty$ then 
$(E_{y_1})^{*(n-1)}=\R^{n-1}$ and, therefore,
\begin{equation}
 \label{gxy.eq.riesz_2}
 \begin{split}
  \int_{E_{x_1}}dx'\int_{E_{y_1}}dy'\,g_{x_1,y_1}&(x'-y')
  \leq 
  \int_{E_{x_1}}dx'\int_{\R^{n-1}}dy'\,g_{x_1,y_1}(x'-y')\\
  &=\mathcal{L}^{n-1}(E_{x_1})\|g_{x_1,y_1}\|_{L^1(\R^{n-1})}
  =\mathcal{L}^{n-1}\big((E_{x_1})^{*(n-1)}\big)\|g_{x_1,y_1}\|_{L^1(\R^{n-1})}\\
  &=\int_{(E_{x_1})^{*(n-1)}}dx'\int_{\R^{n-1}}dy'\,g_{x_1,y_1}(x'-y')\\
  &=\int_{(E_{x_1})^{*(n-1)}}dx'\int_{(E_{y_1})^{*(n-1)}}dy'\,g_{x_1,y_1}(x'-y').
  \end{split}
\end{equation}
In addition, since $g_{x_1,y_1}$ is positive and integrable in $\R^{n-1}$, equality in \eqref{gxy.eq.riesz_2} holds if and only if either 
$\mathcal{L}^{n-1}(E_{x_1})=0$ or 
$E_{y_1}=\R^{n-1}$ up to a set of measure zero. In particular, if $E_{x_1}$ has finite measure, $E_{y_1}$ has infinite measure, and equality in \eqref{gxy.eq.riesz_2} holds, then either $\mathcal{L}^{n-1}(E_{x_1})=0$ or
\begin{equation}\label{eq:case of eq in Riesz cylin_2}
E_{y_1}=(E_{y_1})^{*(n-1)}+c\quad\text{ for every $c\in\R^{n-1}$.}
\end{equation} 
Having these observations in mind, we can now proceed to the proof of the first part of the lemma.

First of all, the assumption $|E\cap((-\pi,\pi)\times\R^{n-1})|<+\infty$ leads to $\mathcal{L}^{n-1}(E_{x_1})<+\infty$ for a.e.\! $x_1\in(-\pi,\pi)$. In particular, 
for a.e.\! $x_1\in(-\pi,\pi)$ and all $y_1\in\R$ with $x_1\neq y_1$, we have
\begin{equation}
\int_{E_{x_1}}dx'\int_{A}dy'\,g_{x_1,y_1}(x'-y')
\leq|E_{x_1}|\|g_{x_1,y_1}\|_{L^1(\R^{n-1})}<+\infty
\end{equation} 
for every measurable set $A\subset\R^{n-1}$. This justifies the forthcoming computations.
Using \eqref{gxy.eq.riesz} and \eqref{gxy.eq.riesz_2} (depending on whether $\mathcal{L}^{n-1}(E_{y_1})<+\infty$ or $\mathcal{L}^{n-1}(E_{y_1})=+\infty$), and recalling \eqref{gy:eq:x1section and rearr interch}, we get
\begin{equation}\label{gy:eq:ineq in cyl rear proof}
\begin{split}
  \PP_s[E]&=\int_{-\pi}^\pi dx_1
  \int_{\R}dy_1
  \int_{E_{x_1}}dx'\int_{\R^{n-1}\setminus E_{y_1}}dy'\,
  g_{x_1,y_1}(x'-y')
  \\
  &=\int_{-\pi}^\pi dx_1
  \int_{\R}dy_1
  \bigg(
  \int_{E_{x_1}}dx'\int_{\R^{n-1}}dy'\,
  g_{x_1,y_1}(x'-y')
  -
  \int_{E_{x_1}}dx'\int_{E_{y_1}}dy'\,
  g_{x_1,y_1}(x'-y')\bigg)
  \\
  &=
  \int_{-\pi}^\pi dx_1
  \int_{\R}dy_1
  \bigg(
  \mathcal{L}^{n-1}(E_{x_1})\int_{\R^{n-1}}dz'\,
  g_{x_1,y_1}(z')-
  \int_{E_{x_1}}dx'\int_{E_{y_1}}dy'\,
  g_{x_1,y_1}(x'-y')\bigg)
  \\
  &\geq 
  \int_{-\pi}^\pi dx_1
  \int_{\R}dy_1
  \bigg(
  \mathcal{L}^{n-1}((E^{*\!\operatorname{cyl}})_{x_1})\int_{\R^{n-1}}dz'\,
  g_{x_1,y_1}(z')\\
   &\hskip100pt-\int_{(E^{*\!\operatorname{cyl}})_{x_1}}\!\!\!\!dx'\int_{(E^{*\!\operatorname{cyl}})_{y_1}}\!\!\!\!dy'\,
  g_{x_1,y_1}(x'-y')\bigg)
  \\
  &=
  \PP_s[E^{*\!\operatorname{cyl}}],
\end{split}
\end{equation}
which proves the first part of the lemma. 

To prove its second part, we need to show that if $E$ has positive finite measure in $(-\pi,\pi)\times\R^{n-1}$ and $\PP_s[E]=\PP_s[E^{*\!\operatorname{cyl}}]<+\infty$, then
there exists $c\in\R^{n-1}$ such that $E_{y_1}=(E_{y_1})^{*(n-1)}+c$ for a.e.\! $y_1\in\R$. 
For this, a key point is to use \eqref{gy:eq:ineq in cyl rear proof}, already proven, as follows. 
Since 
$\mathcal{L}^{n-1}(E_{x_1})
=\mathcal{L}^{n-1}((E^{*\!\operatorname{cyl}})_{x_1})<+\infty$ for a.e.\! $x_1\in(-\pi,\pi)$ and $ \PP_s[E]=\PP_s[E^{*\!\operatorname{cyl}}]<+\infty$,
\eqref{gy:eq:ineq in cyl rear proof} shows that we must have
\begin{equation}
\begin{split}
  \int_{-\pi}^\pi dx_1
  \int_{\R}dy_1
  \bigg(
  \int_{E_{x_1}}&dx'\int_{E_{y_1}}dy'\,
  g_{x_1,y_1}(x'-y')\bigg)\\
  &=
  \int_{-\pi}^\pi dx_1
  \int_{\R}dy_1
  \bigg(
  \int_{(E^{*\!\operatorname{cyl}})_{x_1}}\!\!\!\!dx'\int_{(E^{*\!\operatorname{cyl}})_{y_1}}\!\!\!\!dy'\,
  g_{x_1,y_1}(x'-y')\bigg).
\end{split}
\end{equation}
Using now \eqref{gxy.eq.riesz} in case $\mathcal{L}^{n-1}(E_{y_1})<+\infty$, or \eqref{gxy.eq.riesz_2} in case 
$\mathcal{L}^{n-1}(E_{y_1})=+\infty$, we deduce that for 
a.e.\! $(x_1,y_1)\in (-\pi,\pi)\times\R$ it holds
\begin{equation}
 \label{gxy.eq.riesz.true.eq}
  \int_{E_{x_1}}dx'\int_{E_{y_1}}dy'\,g_{x_1,y_1}(x'-y')
  = 
  \int_{(E_{x_1})^{*(n-1)}}dx'\int_{(E_{y_1})^{*(n-1)}}dy'\,g_{x_1,y_1}(x'-y').
\end{equation}
Now, note that the assumption $0<|E\cap((-\pi,\pi)\times\R^{n-1})|<+\infty$ yields that 
$$
  \Lambda:=\{x_1\in(-\pi,\pi):\,0<\mathcal{L}^{n-1}(E_{x_1})<+\infty\}
$$
has positive measure. 
Recall that we do not assume $E$ to be periodic, and therefore there might exist a set $I\subset(-\infty,-\pi)\cup(\pi,+\infty)$ of positive measure such that  $\mathcal{L}^{n-1}(E_{y_1})=+\infty$ for all $y_1\in I.$  In view of this fact we define 
$$
  \Sigma:=\{y_1\in\R:\,0<\mathcal{L}^{n-1}(E_{y_1})<+\infty\},
$$
which also has positive  
measure because $\Lambda\subset\Sigma$.

Let now $(x_1,y_1)\in\Lambda\times\Sigma$ be such that
\eqref{gxy.eq.riesz.true.eq} holds. Then, by the statement concluding \eqref{eq:case of eq in Riesz cylin} we have 
$$
  E_{x_1}=B_{r(x_1)}^{n-1}+c(x_1,y_1)\quad\text{and}\quad
  E_{y_1}=B_{r(y_1)}^{n-1}+c(x_1,y_1).
$$
Repeating this argument for any other $(\tilde{x}_1,y_1)\in\Lambda\times\Sigma$ we also obtain
$E_{y_1}=B_{r(y_1)}^{n-1}+c(\tilde{x}_1,y_1)$, and hence $c$ is independent of $x_1$. Analogously, the same argument on any other $({x}_1,\tilde y_1)\in\Lambda\times\Sigma$
yields that 
$E_{x_1}=B_{r(x_1)}^{n-1}+c(x_1,\tilde y_1)$. Hence, $c$ is also independent of $y_1$. Thus, there exists $c\in\R^{n-1}$ such that 
$E_{y_1}=B_{r(y_1)}^{n-1}+c$ for a.e.\! $y_1\in \R$ satisfying $0<\mathcal{L}^{n-1}(E_{y_1})<+\infty$.

For this constant $c\in\R^{n-1}$, we now claim that \begin{equation}
 \label{gy:eq:ae y1 in Lambda}
  E_{y_1}=(E_{y_1})^{*(n-1)}+c
\end{equation} 
for a.e.\! $y_1\in\R$. On the one hand, we have already proven \eqref{gy:eq:ae y1 in Lambda} for a.e.\! $y_1\in\Sigma$. On the other hand, if $y_1\in\R\setminus\Sigma$, then we must have that either $\mathcal{L}^{n-1}(E_{y_1})=0$ or $\mathcal{L}^{n-1}(E_{y_1})=+\infty$. But if $\mathcal{L}^{n-1}(E_{y_1})=0$, then $E_{y_1}=(E_{y_1})^{*(n-1)}+d$ for every $d\in\R^{n-1}$ up to a set of $\mathcal{L}^{n-1}$-measure zero, since 
$(E_{y_1})^{*(n-1)}$ is the empty set in 
$\R^{n-1}$. Thus, \eqref{gy:eq:ae y1 in Lambda} holds in this case too.
It only remains to show \eqref{gy:eq:ae y1 in Lambda} 
for a.e.\! $y_1\in\R$ such that $\mathcal{L}^{n-1}(E_{y_1})=+\infty$. Recall that \eqref{gxy.eq.riesz.true.eq} holds for 
a.e.\! $(x_1,y_1)\in (-\pi,\pi)\times\R$. Hence, since $|\Lambda|>0$, for 
a.e.\! $y_1\in\R$ with $\mathcal{L}^{n-1}(E_{y_1})=+\infty$ there exists a point $x_1\in\Lambda$ such that \eqref{gxy.eq.riesz.true.eq} holds. Then the statement concluding \eqref{eq:case of eq in Riesz cylin_2} yields that 
$E_{y_1}=(E_{y_1})^{*(n-1)}+d$
for every $d\in\R^{n-1}$. Thus, \eqref{gy:eq:ae y1 in Lambda} also follows in this last case.

In conclusion, there exists $c\in\R^{n-1}$ such that $E_{y_1}=(E_{y_1})^{*(n-1)}+c$ for a.e.\! $y_1\in\R$.
Equivalently, $E=E^{*\!\operatorname{cyl}}+(0,c)$ up to a set of measure zero.  
\end{proof}


We now address the proof of Theorem \ref{thm:intro:symmetry of Delaunay}. Let us first recall some definitions to facilitate the reading. For a measurable set $A\subset\R^n$ we denote its Steiner symmetrization with respect to 
$\{x_1=0\}$ by $A^*.$ If $A\subset \R$ is measurable and of finite measure, then $A^*$ is an interval centered at $0$ with the same measure as $A.$
The Steiner symmetrization
is well defined whenever the $1$-dimensional set $A\cap\{x'=c\}$ has finite length for all $c\in\R^{n-1}.$ Thus, the periodic rearrangement
$$
  E^{*\!\eper}:=\bigcup_{k\in\Z}\left[\left(E\cap \{-\pi<x_1<\pi\}\right)^*+2k\pi e_1\right]
$$
is well defined for every measurable set $E\subset\R^n$ which is $2\pi$-periodic in the $x_1$ direction. 

By periodicity, note that $ E^{*\!\eper}= 
(E+(c,0))
^{*\!\eper}$ for all $c\in\R.$

We denote the symmetric decreasing rearrangement of a real-valued function $f$ on $\R$ by $f^*$. It is defined by the following properties: $f^*$ is nonnegative, equimeasurable with $|f|$, and its superlevel sets are open intervals centered at $0$. Two explicit expressions for $f^{*}$ are
\begin{equation}
 \label{eq:gy:layer cake for general rearrange}
  f^*(x)=\int_0^{+\infty}dt\,\chi_{\{|f|>t\}^*}(x)=\sup\big\{t:\,x\in\{|f|>t\}^*\big\}.
\end{equation}
For a $2\pi$-periodic function $u:\R\to\R$, the periodic symmetric deceasing rearrangement of $u$, denoted by $u^{*\!\eper}$, is the unique $2\pi$-periodic function such that 
$u^{*\!\eper}\chi_{(-\pi,\pi)}=(u\chi_{(-\pi,\pi)})^*$. 

It follows from these definitions that if $E$ is as in \eqref{Del frac min E u}, then
\begin{equation}\label{eq:gy:E asteper and u asteper}
  E^{*\!\eper}=\{x\in\R^n\,:\,|x'|<u^{*\!\eper}(x_1)\}.
\end{equation}

The main ingredient of our proof of \Cref{thm:intro:symmetry of Delaunay} will be the following Riesz rearrangement inequality on the circle.

\begin{theorem}$($\cite[Theorem 2]{Baernstein Taylor}, \cite[Theorem 1]{FriedLutt}, \cite[Theorem 2]{Burchard Hajaiej}, \cite[Theorem 7.3]{Baernstein book}$)$
\label{thm:gy:sharp Friedberg Luttinger}
Let $f,h,g:\R\to \R$ be three nonnegative $2\pi$-periodic measurable functions. Assume that $g$ is even, as well as nonincreasing in $(0,\pi)$. 

Then, 
\begin{equation}\label{eq:thm FriedLutt inequality}
    \int_{-\pi}^{\pi}dx\int_{-\pi}^{\pi} dy\,f(x)g(x-y)h(y)
    \leq 
    \int_{-\pi}^{\pi}dx\int_{-\pi}^{\pi}dy\, f^{*}(x)g(x-y)h^{*}(y).
\end{equation}
In addition, if $g$ is decreasing in $(0,\pi)$ 
and the left-hand side of \eqref{eq:thm FriedLutt inequality} is finite,
then
equality holds in \eqref{eq:thm FriedLutt inequality} if and only if at least one of the following conditions holds: 
\begin{itemize}
\item[$(i)$] either $f$ or $h$ is constant almost everywhere.

\item[$(ii)$] there exists $z\in \R$ such that $f(x)=f^{*\!\eper}(x+z)$ and $h(x)=h^{*\!\eper}(x+z)$ for almost every $x\in\R$.
\end{itemize}
\end{theorem}

The inequality in \Cref{thm:gy:sharp Friedberg Luttinger} was first discovered, independently, in \cite{Baernstein Taylor} and \cite{FriedLutt}. Both references contain more general inequalities: \cite{Baernstein Taylor} deals with the sphere $\mathbb{S}^n$, whereas \cite{FriedLutt} deals with a product of more than three functions. 
The result \cite[Theorem 7.3]{Baernstein book} is also more general than Theorem \ref{thm:gy:sharp Friedberg Luttinger}, as it deals with a Riesz rearrangement inequality on the sphere $\mathbb{S}^n$ and not only on the circle. Moreover, \cite{Baernstein book} treats more general functions of $f(x)$ and $h(y)$ than simply the product $f(x)h(y).$ 
The inequality \eqref{eq:thm FriedLutt inequality} can also be found in \cite{Baernstein Circle} in a more general form where $g$ is also rearranged.

Instead, the statement in \Cref{thm:gy:sharp Friedberg Luttinger} concerning equality in \eqref{eq:thm FriedLutt inequality} follows from Burchard and Hajaiej \cite[Theorem 2]{Burchard Hajaiej}, who treated the case of equality in $\mathbb{S}^n$ for the first time. For this result, we also cite \cite[Theorem 7.3]{Baernstein book} since, being less general than \cite{Burchard Hajaiej}, fits precisely with our setting.

We find reference \cite{Baernstein book} to be the simplest one for looking up all statements of \Cref{thm:gy:sharp Friedberg Luttinger}.

\begin{proof}[Proof of Theorem \ref{thm:intro:symmetry of Delaunay}.]
Let us first show that $\PP_{s}[E^{*\!\eper}]\leq \PP_{s}[E].$
In contrast to the proof of Lemma \ref{gy:lemma:cylindrical rearr}, we now write $\PP_{s}[E]$ as
$$
    \PP_{s}[E]=\int_{\R^{n-1}}dx'\int_{E_{x'}\cap(-\pi,\pi)}dx_1
    \int_{\R^{n-1}}dy'\int_{\R\setminus E_{y'}}dy_1\,\frac{1}{|x-y|^{n+s}},
$$
where we have set $E_{x'}:=\{x_1\in \R:\,(x_1,x')\in E\}$ (in analogy with the sections $E_{x_1}$ defined above).
Writing $\PP_{s}[E^{*\!\eper}]$ in the same way (simply replacing $E$ by $E^{*\!\eper}$) and interchanging the order of integration $dx_1\,dy'$, we see that it is sufficient to prove that
\begin{equation}
\label{eq:sections x2 and y2}
    \int_{E_{x'}\cap(-\pi,\pi)}dx_1
    \int_{\R\setminus E_{y'}}
    \frac{dy_1 }{|x-y|^{n+s}}
    \geq
      \int_{(E^{*\!\eper})_{x'}\cap(-\pi,\pi)}dx_1\int_{\R\setminus(E^{*\!\eper})_{y'}} \frac{dy_1}{|x-y|^{n+s}}
\end{equation}
for almost every $(x',y')\in \R^{n-1}\times\R^{n-1}$.

For this, set
$$
    a=|x'-y'|\quad\text{ and }\quad \gamma=\frac{n+s}{2} \quad\text{ (and thus } |x-y|^{n+s}=\big((x_1-y_1)^2+a^2\big)^{\gamma}\text{)},
$$
observe that the following integrals are all finite when $a\neq 0$, and that the condition $a\neq 0$ holds for almost every $(x',y')\in \R^{n-1}\times\R^{n-1}$. We have 
$$
   \int_{\R\setminus E_{y'}} \frac{dy_1}{|x-y|^{n+s}}
   =
   \int_{\R} \frac{dy_1}{|x-y|^{n+s}}
   -\int_{E_{y'}} \frac{dy_1}{|x-y|^{n+s}}
$$
and
$$
   \int_{\R\setminus(E^{*\!\eper})_{y'}} \frac{dy_1}{|x-y|^{n+s}}
   =
   \int_{\R} \frac{dy_1}{|x-y|^{n+s}}
   -\int_{(E^{*\!\eper})_{y'}} \frac{dy_1}{|x-y|^{n+s}}.
$$
From the equality $|E_{x'}\cap (-\pi,\pi)|=|(E^{*\!\eper})_{x'}\cap (-\pi,\pi)|,$ we infer that
$$
    \int_{E_{x'}\cap(-\pi,\pi)}dx_1\int_{\R} \frac{dy_1}{|x-y|^{n+s}}
    =
    \int_{(E^{*\!\eper})_{x'}\cap(-\pi,\pi)}dx_1\int_{\R}\frac{dy_1}{|x-y|^{n+s}}.
$$
Therefore, to prove  \eqref{eq:sections x2 and y2} it is enough to show that
\begin{equation}\label{eq:sections x2 and y2 without complement}
    \int_{E_{x'}\cap(-\pi,\pi)}dx_1\!\int_{E_{y'}}\!\frac{dy_1}{((x_1-y_1)^2+a^2)^{\gamma}}
    \leq 
    \int_{(E^{*\!\eper})_{x'}\cap(-\pi,\pi)}dx_1\!\int_{(E^{*\!\eper})_{y'}}\! \frac{dy_1}{((x_1-y_1)^2+a^2)^{\gamma}}
\end{equation}
for all $a\neq 0$.

In order to do this, recall that $(E^{*\!\eper})_{y'}=(E_{y'})^{*\!\eper}$ and that  both $E_{y'}$ and $(E^{*\!\eper})_{y'}$ are $2\pi$-periodic. Taking these facts in consideration, \eqref{eq:sections x2 and y2 without complement} will follow if we show that, for every two measurable subsets $A$ and $B$ of $(-\pi,\pi)$, it holds
\begin{equation}
\begin{split}
    \int_{A}dx_1\int_{\bigcup_{k\in\Z}(B+2k\pi)} 
    &\frac{dy_1}{((x_1-y_1)^2+a^2)^{\gamma}}
    \leq
    \int_{A^{*}}dx_1\int_{\bigcup_{k\in\Z}(B^{*}+2k\pi)}
    \frac{dy_1}{((x_1-y_1)^2+a^2)^{\gamma}}.
\end{split}
\end{equation}
By the change of variables $\overline{y}_1=y_1-2k\pi$, this can be equivalently written as
\begin{equation}\label{int_Aint_B g}
    \int_A dx_1\int_B d\overline y_1\,g(x_1-\overline y_1)\leq 
    \int_{A^{*}}dx_1\int_{B^{*}}d\overline y_1\,g(x_1-\overline y_1)
    \text{ for all subsets $A$ and $B$ of $(-\pi,\pi)$},
\end{equation}
where we have set
\begin{equation}\label{def.g.rear}
    g(z):=\sum_{k\in\Z}\frac{1}{((z+2k\pi)^2+a^2)^{\gamma}}.
\end{equation}
Clearly, the function $g$ is positive, $2\pi$-periodic, and even. We will next show that $g$ is decreasing in $(0,\pi)$, and hence \eqref{int_Aint_B g} will follow from Theorem \ref{thm:gy:sharp Friedberg Luttinger} applied to the characteristic functions $f=\chi_A$ and $h=\chi_B$.

In order to show that $g$ is decreasing in $(0,\pi)$, we introduce the Laplace transform of the function $t\mapsto t^{\gamma-1}$, which amounts to the equality
\begin{equation}\label{gamma_x-s}
  w^{-\gamma}=\frac{1}{\Gamma(\gamma)}\int_0^{+\infty}dt\,
  t^{\gamma-1}e^{-wt}\quad\text{for all $w>0$,}
\end{equation}
where 
$\Gamma(\gamma):=\int_0^{+\infty}dr\,{r}^{\gamma-1}e^{-r}$
is the Gamma function.
The identity \eqref{gamma_x-s} follows simply by a change of variables $r=wt$ in the definition of $\Gamma(\gamma)$.
Using \eqref{gamma_x-s} in the definition of $g$ yields
\begin{equation}\label{g.laplace.heat}
 g(z)=\frac{1}{\Gamma(\gamma)}\sum_{k\in \Z}\int_0^{+\infty}dt\,
 t^{\gamma-1}e^{-((z+2k\pi)^2+a^2)t}
 =\frac{1}{\Gamma(\gamma)}\int_0^{+\infty}dt\,t^{\gamma-1}e^{-a^2t}G_t(z),
\end{equation}
where
$$
  G_t(z):=\sum_{k\in \Z}e^{-(z+2k\pi)^2t}.
$$
We will now show that $G_t$ is decreasing in $(0,\pi)$ for every $t>0$, yielding, in view of \eqref{g.laplace.heat}, that $g$ is also decreasing in $(0,\pi).$ The claim on the strict monotonicity of $G_t$ follows directly from the fact that the fundamental solution of the heat equation in $(-\pi,\pi)$ with periodic boundary conditions (or, in other words, of the heat equation on the circle), which is 
$$
  \Lambda_t(z):=\frac{1}{\sqrt{4\pi t}}\sum_{k\in\mathbb{Z}}
  e^{-\frac{(z+2k\pi)^2}{4t}},
$$
is decreasing in $(0,\pi)$ for all $t>0$. For the sake of completeness, in \Cref{thm:appendix:decreasing} we give a simple proof of this well-known fact, based solely on maximum principles for the heat equation.


We now prove that if $\PP_{s}[E^{*\!\eper}]=\PP_{s}[E]<+\infty$, then  $E=E^{*\!\eper}+ce_1$ for some $c\in\R$, up to a set of measure zero. Indeed, in view of \eqref{eq:sections x2 and y2} and the equalities following it, if $\PP_{s}[E^{*\!\eper}]=\PP_{s}[E]$ then there must be equality in \eqref{eq:sections x2 and y2 without complement} for almost every $(x',y')\in \R^{n-1}\times \R^{n-1}$. Equivalently, using the formulation of \eqref{int_Aint_B g},
\begin{equation}\label{int_Aint_B g.eq.case}
\begin{split}
    \int_{E_{x'}\cap(-\pi,\pi)}\!\!\!dx_1
    \int_{E_{y'}\cap(-\pi,\pi)}\!\!\!dy_1\, g(x_1-y_1)
    = 
    \int_{(E^{*\!\eper})_{x'}\cap(-\pi,\pi)}\!\!\!dx_1
    \int_{(E^{*\!\eper})_{y'}\cap(-\pi,\pi)}\!\!\!dy_1\,g(x_1-y_1)
    \end{split}
\end{equation}
for almost every $(x',y')\in \R^{n-1}\times \R^{n-1}$. 

Now, let $J\subset\R^{n-1}$ be the set of points $x'\in \R^{n-1}$ for which $E_{x'}\cap(-\pi,\pi)$ is neither $(-\pi,\pi)$ nor $\varnothing$ (up to zero measure sets). We may assume that $J$ has positive $(n-1)$-dimensional Lebesgue measure ---if $J$ has measure zero then $E=E^{*\!\eper}$ up to a set of measure zero in $\R^n$, and we would be done. 
Then, $J\times J$ has positive measure in $\R^{n-1}\times \R^{n-1}$, and \eqref{int_Aint_B g.eq.case} holds for almost every $(x',y')\in J\times J$. Note that, by definition of $J$, we don't fall in case $(i)$ of Theorem \ref{thm:gy:sharp Friedberg Luttinger} when using \eqref{int_Aint_B g.eq.case} on points $x'$ and $y'$ in $J$. Hence, by \eqref{int_Aint_B g.eq.case} and \Cref{thm:gy:sharp Friedberg Luttinger} $(ii)$ we deduce that, for almost every $(x',y')\in J\times J$, we have
$$
  E_{x'}=(E^{*\!\eper})_{x'}+c(x',y')\quad\text{ and }\quad E_{y'}=(E^{*\!\eper})_{y'}+c(x',y')
$$
for some constant $c(x',y')\in\R$ depending, a priori, on $x'$ and $y'$.
This easily yields that $c(x',y')=c$ is actually independent of $x'$ and $y'$. Indeed, considering two different values $x'\in J$ and $\overline{x}'\in J$ for the same $y'$, we have
\begin{equation}\label{gy:eq:cxy indep}
  (E^{*\!\eper})_{y'}+c(x',y')=E_{y'}=(E^{*\!\eper})_{y'}+c(\overline{x}',y').
\end{equation}
Since $y'\in J$, $(E^{*\!\eper})_{y'}$ is a $2\pi$-periodic set for which $(E^{*\!\eper})_{y'}\cap(-\pi,\pi)$ is a nonempty open
interval different from $(-\pi,\pi)$. Therefore \eqref{gy:eq:cxy indep} yields that $c(x',y')=c(\overline{x}',y')$ modulo $2\pi$. By symmetry, $c(x',y')$ is also independent of $y'$. Therefore, 
$E_{y'}=(E^{*\!\eper})_{y'}+c$ for some $c\in\R$ and all $y'\in J$. By the definition of $J$, we conclude that $E=E^{*\!\eper}+ce_1$ up to a set of measure zero.
\end{proof}


\begin{remark}{\em
From the last part of the proof of \Cref{thm:intro:symmetry of Delaunay}, in particular in what regards \eqref{gamma_x-s}, we see that  \Cref{thm:intro:symmetry of Delaunay} easily extends to any kernel which is the Laplace transform of a nonnegative function. More precisely, consider a periodic nonlocal perimeter functional of the form
$$
  \PP_{K}[E]:=\int_{E\cap\{-\pi<x_1<\pi\}}dx\int_{\R^n\setminus E}dy\,K(|x-y|^2),
$$
where $K:(0,+\infty)\to(0,+\infty)$ is a kernel given by 
$$ K(w):=\int_0^{+\infty}dt\,\kappa(t)e^{-wt},\quad w>0,$$
for some measurable function $\kappa: (0,+\infty)\to[0,+\infty)$ which is not identically zero. Then, Theorem \ref{thm:intro:symmetry of Delaunay} holds also for $\PP_{K}$. 

Indeed, arguing as in the  proof of Theorem \ref{thm:intro:symmetry of Delaunay} up to \eqref{int_Aint_B g}, one sees that it is enough to show that 
\begin{equation}
    \int_A dx_1\int_B dy_1\, g_K(x_1-y_1)\leq 
    \int_{A^{*}}dx_1 \int_{B^{*}}dy_1\,g_K(x_1-y_1)
    \quad\text{for all subsets $A$ and $B$ of $(-\pi,\pi)$},
\end{equation}
where 
\begin{equation}
    g_K(z):=\sum_{k\in\Z}K((z+2k\pi)^2+a^2),\quad z\in\R.
\end{equation}
But
\begin{equation}
 g_K(z)=\sum_{k\in \Z}\int_0^{+\infty}dt\,
 \kappa(t)e^{-((z+2k\pi)^2+a^2)t}
 =\int_0^{+\infty}dt\,\kappa(t)e^{-a^2t}G_t(z)
\end{equation}
is a decreasing function in $(0,\pi)$ by \Cref{thm:appendix:decreasing} and by  the fact that $\kappa\geq0$ is not identically zero in 
$(0,+\infty)$. Hence, Theorem \ref{thm:gy:sharp Friedberg Luttinger} still applies in this framework and, thus, we deduce that $\PP_K$ is nonincreasing under periodic decreasing rearrangement. Moreover, as before, if
$\PP_K[E^{*\!\eper}]= \PP_K[E]$ then $E=E^{*\!\eper}+ce_1$ for some $c\in\R$, up to a set of measure zero.
}\end{remark}

\appendix


\section{Periodic minimizers in $\R^2$}\label{ss CFSW n2}
\label{section:periodic min in R2}

In this appendix we present, when $n=2$, an expression for the functional $\PP_s$  which is simpler than that of \Cref{scan lemma 2} in $\R^n$. As mentioned in the Introduction, this expression allowed us to prove the existence of minimizer when $n=2$ in a simple way. We later found the existence proof of Theorem \ref{scan thm1} in all dimensions.

The formula for the periodic version of the fractional perimeter when $n=2$ reads as follows.
If $E=\{(x_1,x_2)\in\R^2:\,|x_2|<u(x_1)\}$, for some nonnegative $2\pi$-periodic function $u$, then $\PP_s[E]$ can be written in terms of $u$ as
\begin{equation}\label{scan n2}
\PP_s[E]=\int_{-\pi}^\pi dx_1\!
\int_{-\infty}^{+\infty}dy_1\,\frac{2}{|x_1-y_1|^{s}}
\Big\{G\Big(\frac{u(x_1)-u(y_1)}{|x_1-y_1|}\Big)
+H\Big(\frac{u(x_1)+u(y_1)}{|x_1-y_1|}\Big)\Big\},
\end{equation}
where
\begin{equation}
G(t):=\int_0^t d\tau\, (t-\tau)(1+\tau^2)^{-\frac{2+s}{2}}
\quad\text{and}\quad
H(t) := G'(+\infty)t-G(t)\quad\text{for $t\in\R$.}
\end{equation}
Note that $G(0)=0$,
\begin{equation}\label{F}
0<G'(t)=\int_0^td\tau\, (1+\tau^2)^{-\frac{2+s}{2}}
<G'(+\infty) < +\infty,
\end{equation}
and
$G''(t)=(1+t^2)^{-\frac{2+s}{2}}>0$ for all $t\in\R$.
Therefore, $G$ is a nonnegative even function which is increasing in
$(0,+\infty)$ and strictly convex on $\R$. Furthermore, it behaves as a quadratic function near the origin and approaches a linear one as $t\to\pm\infty$.
On the other hand, $H(0)=0$, $H'=G'(+\infty)-G'>0$, and $H''=-G''<0$ in $[0,+\infty).$
Thus, $H$ is positive, increasing, and strictly concave in $[0,+\infty)$. Finally, by Fubini's theorem,
\begin{equation}
\begin{split}
H(+\infty)=\int_0^{+\infty}dt\,H'(t)&=\int_0^{+\infty}dt\,(G'(+\infty)-G'(t))
=\int_0^{+\infty}dt\int_t^{+\infty}d\tau\, (1+\tau^2)^{-\frac{2+s}{2}}
\\
&=\int_0^{+\infty}d\tau\,\tau(1+\tau^2)^{-\frac{2+s}{2}}<+\infty.
\end{split}
\end{equation}

For the sake of completeness we will give the proof of  \eqref{scan n2}, even though this expression was not used to prove the results of this paper. But before, let us explain how to use \eqref{scan n2} to give a simple proof of the existence of a constrained minimizer of $\PP_s$ when $n=2$. The argument goes as follows. Since $G(t)+C_1\geq C_2|t|$ for some constants $C_1,C_2>0$, the term  in  \eqref{scan n2} involving
 $G$  bounds the $W^{s,1}(-\pi,\pi)$ norm of $u$ up to an additive constant.
In this way, since $H\geq 0$, a minimizing sequence for $\PP_s$ is a bounded sequence in $W^{s,1}(-\pi,\pi)$. By compactness, this yields the existence of the desired constrained minimizer. We refer to \cite{Alvinya} for the details on this argument; see also the end of Part~1 in our proof of Theorem \ref{scan thm1}.

\begin{lemma}\label{lemma heyprop}
Let $u:\R\to\R$ be nonnegative, $2\pi$-periodic, and such that  $u\in W^{s,1}(-\pi-\epsilon,\pi+\epsilon)$ for some (and hence for all)
$\epsilon>0$. Let 
$E=\{(x_1,x_2)\in\R^2:\,|x_2|<u(x_1)\}$. Then, 
\eqref{scan n2} holds.
\end{lemma}
\begin{proof}
We prove \eqref{scan n2} using \eqref{p2}. Let us first rewrite $\phi$ in \eqref{p3} for $n=2$ in a more convenient way by introducing the function $G'$. For $p\geq 0$ and $q\geq0$ we have (using the change of variables $z=\tau+w$ twice)
\begin{equation}
\begin{split}
\phi(p,q)
&=\int_{-p}^pdw\,\Big\{\int_{-\infty}^{-q}dz\,(1+|z-w|^2)^{-\frac{2+s}{2}}
+\int^{+\infty}_{q}dz\,(1+|z-w|^2)^{-\frac{2+s}{2}}\Big\}\\
&=\int_{-p}^pdw\,\Big\{G'(-q-w)-G'(-\infty)+G'(+\infty)-G'(q-w)\Big\}.
\end{split}
\end{equation}
Note that $G'$ is odd, $0<-G'(-\infty)=G'(+\infty)<+\infty$, and $G'(-q-\cdot)-G'(q-\cdot)$ is an even function. Therefore,
\begin{equation}\label{p2 aux1}
\begin{split}
\phi(p,q)&=2\int_{0}^pdw\,\Big\{G'(-q-w)-G'(q-w)+2G'(+\infty)\Big\}\\
&=4G'(+\infty)p+2\int_{0}^pdw\,G'(-q-w)-2\int_{0}^pdw\,G'(q-w)\\
&=4G'(+\infty)p-2\int_{-q}^{-q-p}dt\,G'(t)+2\int_{q}^{q-p}dt\,G'(t)\\
&=2\big(2G'(+\infty)p+G(p-q)-G(p+q)\big),
\end{split}
\end{equation}
where we also used in the last equality that $G$ is even. Now, from \eqref{p2} and \eqref{p2 aux1} we deduce that
\begin{equation}\label{hey1aa}
\begin{split}
\PP_s[E]& =2\int_{-\pi}^\pi dx_1
\int_{-\infty}^{+\infty}dy_1\,|x_1-y_1|^{-s}
\\
&\phantom{aaaaaaaa}\times\Big\{G\Big(\frac{u(x_1)-u(y_1)}{|x_1-y_1|}\Big)
-G\Big(\frac{u(x_1)+u(y_1)}{|x_1-y_1|}\Big)+2G'(+\infty)\,
\frac{u(x_1)}{|x_1-y_1|}\Big\}.
\end{split}
\end{equation}

Our goal is to show that in the term involving
$2G'(+\infty)u(x_1)|x_1-y_1|^{-1}$ we can replace $u(x_1)$ by $u(y_1)$ and the identity \eqref{hey1aa} remains unchanged. Once this is proved, we get that in this last term one can replace $2u(x_1)$ by $u(x_1)+u(y_1)$. Hence, using that $H(t) = G'(+\infty)t-G(t)$, this gives \eqref{scan n2}. 

Therefore, it only remains to show that we can replace $u(x_1)$ by $u(y_1)$ in the last term on the right-hand side of \eqref{hey1aa}. 
This follows immediately from \Cref{gy:lemma:interch integ periodic} (used with $m=l=1$) applied to the integral
\begin{equation}
     \int_{-\pi}^\pi dx\!\int_{-\infty}^{+\infty}dy\,\frac{f(x,y)}{|x-y|^{s}},
\end{equation}
where $$f(x,y):=-G\Big(\frac{u(x)+u(y)}{|x-y|}\Big)
     +2G'(+\infty)\,\frac{u(x)}{|x-y|}.$$
To check the hypothesis of the lemma, we need to verify that     
\begin{equation}
  \label{gy:eq:f integral bound}
     \int_{-\pi}^\pi dx\!\int_{-\infty}^{+\infty}dy\,
     \frac{|f(x,y)|}{|x-y|^{s}}<+\infty.
\end{equation} 
For this purpose we write $f(x,y)=A+B,$ where
$$
  A:=-G\Big(\frac{u(x)+u(y)}{|x-y|}\Big)
  +G\Big(\frac{2u(x)}{|x-y|}\Big)
$$
and
$$
  B:=-G\Big(\frac{2u(x)}{|x-y|}\Big)+2G'(+\infty)\,
\frac{u(x)}{|x-y|}=H\Big(\frac{2u(x)}{|x-y|}\Big).
$$
Using that $u$ is nonnegative, $H$ is nonnegative and increasing in $(0,+\infty)$, and $H(+\infty)<+\infty$, we have 
\begin{equation}
  \label{gy:eq:f integral bound.aux1}
0\leq B\leq H(+\infty)<+\infty.
\end{equation} 
Concerning the quantity $A$, observe that
\begin{equation}
  \label{gy:eq:f integral bound.aux2}
  |A|=\left|  \int_{\frac{u(x)+u(y)}{|x-y|}}^{\frac{2u(x)}{|x-y|}}
  dt\,G'(t)\right|\leq \sup_{t\in\R}|G'(t)|\,\frac{|u(x)-u(y)|}{|x-y|}=G'(+\infty)\frac{|u(x)-u(y)|}{|x-y|}.
\end{equation} 

Recall that we are assuming $u\in W^{s,1}(-\pi-\epsilon,\pi+\epsilon)$ for some $\epsilon>0$. According to this, to prove \eqref{gy:eq:f integral bound} we split the integral and we use \eqref{gy:eq:f integral bound.aux1} and \eqref{gy:eq:f integral bound.aux2} to get
\begin{equation}\label{gy:eq:x minus y geq 3pi.aux3}
\begin{split}
    &\int_{-\pi}^\pi dx\!\int_{-\infty}^{+\infty}dy\,\frac{|f(x,y)|}{|x-y|^{s}}
    =
    \int_{-\pi}^\pi dx\!\int_{|x-y|\leq\epsilon}dy\,\frac{|A|+|B|}{|x-y|^{s}}
    +
    \int_{-\pi}^\pi dx\!\int_{|x-y|>\epsilon}dy\,\frac{|f(x,y)|}{|x-y|^{s}}
    \\
    &\leq 
    G'(+\infty)\|u\|_{W^{s,1}(-\pi-\epsilon,\pi+\epsilon)}
    +H(+\infty)\int_{-\pi}^\pi dx\!\int_{-\pi-\epsilon}^{\pi+\epsilon}
    \frac{dy}{|x-y|^{s}}
    + \int_{-\pi}^\pi dx\!\int_{|x-y|>\epsilon}dy\,\frac{|f(x,y)|}{|x-y|^{s}}.
\end{split}
\end{equation}
The first two terms on the right-hand side of \eqref{gy:eq:x minus y geq 3pi.aux3} are finite.
To estimate the third one, we use that $0\leq G(t)\leq G'(+\infty)t$ for all $t\geq 0,$ and hence 
\begin{equation}\label{gy:eq:x minus y geq 3pi.aux4}
 \begin{split}
    \int_{-\pi}^\pi dx\!\int_{|x-y|>\epsilon}dy\,\frac{|f(x,y)|}{|x-y|^{s}}
    \leq &
    \,G'(+\infty)\int_{-\pi}^\pi dx\!\int_{|x-y|>\epsilon}dy\,\frac{3u(x)}{|x-y|^{1+s}}
    \\
    &+G'(+\infty)
    \int_{-\pi}^\pi dx\!\int_{|x-y|>\epsilon}dy\,\frac{u(y)}{|x-y|^{1+s}}.
\end{split}
\end{equation}
The first double integral in the right-hand side of \eqref{gy:eq:x minus y geq 3pi.aux4} is finite, since $s>0$ and $u\in W^{s,1}(-\pi-\epsilon,\pi+\epsilon)\subset L^1(-\pi,\pi)$.
For the second double integral, since $u$ is $2\pi$-periodic and 
nonnegative, using \Cref{gy:lemma:interch integ periodic} we deduce
\begin{equation}
\int_{-\pi}^\pi dx\!\int_{|x-y|>\epsilon}dy\,\frac{u(y)}{|x-y|^{1+s}}
=\int_{-\pi}^\pi dx\!\int_{\R}dy\,\frac{\chi_{\{|x-y|>\epsilon\}}}{|x-y|^{1+s}}\,u(y)
=\int_{-\pi}^{\pi}d\overline y\,u(\overline y)\int_{|\overline x-\overline y|>\epsilon} \frac{d\overline x}{|\overline x-\overline y|^{1+s}}.
\end{equation}
As before, this is finite since $s>0$ and 
$u\in L^1(-\pi,\pi)$.
Hence  \eqref{gy:eq:f integral bound} follows from \eqref{gy:eq:x minus y geq 3pi.aux3} and the finiteness of the integrals in  \eqref{gy:eq:x minus y geq 3pi.aux4}.
\end{proof}



\section{Strict monotonicity of the periodic heat kernel}
\label{Appendix:Heat Kernel}

In this section we address the proof of an important  ingredient in the proof of \Cref{thm:intro:symmetry of Delaunay}, namely, that the function defined in \eqref{def.g.rear} is decreasing in $(0,\pi)$. As shown there, this property follows from the fact that the periodic heat kernel is decreasing with respect to the spatial variable in $(0,\pi)$. This  is a well-known fact \cite{ACMM18,Andersson,Chavel,CheYa,NSS}, whose proofs are commented next. Here, for the sake of completeness,  we provide a simple proof based on maximum principles for the heat equation. 

Prior to this, let us briefly comment on the proofs given in the above-mentioned references. The one in \cite{ACMM18} follows exactly the same lines as the one presented here, but it is carried out within the context of hypersurfaces of revolution. The one in \cite{CheYa} is also based on the strong maximum principle for the heat equation, but it is more technically involved than \cite{ACMM18}, since a more general class of manifolds is considered there. The proof in 
\cite{Andersson} is based on an explicit computation. The one in \cite{Chavel} follows by expressing the derivative of the kernel in terms of the heat kernel in smaller dimensions and using an induction argument; see \cite[\S6.3]{Chavel}. Finally, the one in \cite{NSS} is a straightforward application of the sharp estimates for the heat kernel on the sphere $\mathbb{S}^n$ found there.

The fundamental solution of the heat equation in $(-\pi,\pi)$ with periodic boundary conditions is given by 
\begin{equation}\label{eq:def of Gammat}
   \Gamma(z,t):=\frac{1}{\sqrt{4\pi t}}\sum_{k\in \Z}e^{-\frac{(z+2k\pi)^2}{4t}} \quad\text{for $z\in\R$ and $t>0$}.
\end{equation}
Given a $2\pi$-periodic initial data $g,$ the function 
\begin{equation}
u(x,t)=\int_{-\pi}^{\pi}dy\,\Gamma(x-y,t)g(y)
\end{equation} 
satisfies 
$\partial_t u-\partial_{xx}u=0$ in $\R\times (0,+\infty)$, that $u(\cdot,t)$ is $2\pi$-periodic for all $t>0$, and that 
$u(\cdot,0)=g$ in $\R$.
This follows immediately from the properties of the standard (nonperiodic) heat kernel and the fact that
\begin{equation}
\begin{split}
  \sum_{k\in \Z}\int_{-\pi}^{\pi}dy\,e^{-\frac{(x-y+2k\pi)^2}{4t}}g(y)
  &=
  \sum_{k\in \Z}\int_{-\pi}^{\pi}dy\,e^{-\frac{(x-y+2k\pi)^2}{4t}}g(y-2k\pi)\\
  &=
  \sum_{k\in \Z}\int_{-\pi-2k\pi}^{\pi-2k\pi}d\overline y\,e^{-\frac{(x-\overline y)^2}{4t}}g(\overline y)
  =
  \int_{\R}d \overline y\,e^{-\frac{(x-\overline y)^2}{4t}}g(\overline y).
\end{split}
\end{equation}

We now state and prove the result on the periodic heat kernel that we used in the proof of \Cref{thm:intro:symmetry of Delaunay}.

\begin{theorem}
\label{thm:appendix:decreasing}
For every $t>0$, the function
$$
   z\mapsto \Gamma(z,t)\quad 
   \text{is decreasing in }(0,\pi).
$$
\end{theorem}

\begin{proof}
Let us first show that $\Gamma(\cdot,t)$ is nonincreasing in $(0,\pi)$. For this, we take functions $g_\epsilon\in C_c^{\infty}(-\epsilon,\epsilon)$, where $0<\epsilon<\pi/2$, approximating the Dirac delta as $\epsilon\downarrow0$, and we consider their $2\pi$-periodic extensions from $[-\pi,\pi]$ to $\R$ (which we also denote by $g_\epsilon$). In particular, we assume that
$$
  \int_{-\epsilon}^{\epsilon} g_\epsilon=1\quad \text{for all $0<\epsilon<\pi/2$}.
$$
We additionally assume that the $g_\epsilon$ are nonnegative, even, and nonincreasing in $(0,\pi)$.

Let $u_{\epsilon}$ be the solution to 
$$
  \left\{\begin{array}{ll}
           \partial_t u_{\epsilon}-\partial_{xx}u_{\epsilon}=0 &\text{ in }\R\times (0,+\infty),
           \smallskip \\
           u_{\epsilon}(\cdot,0)=g_\epsilon & \text{ in $\R$}.
         \end{array}\right.
$$
Thus, $u_{\epsilon}(x,t)=\int_{-\pi}^{\pi}dy\,\Gamma(x-y,t)g_{\epsilon}(y)$ and $u_{\epsilon}(\cdot,t)$ is $2\pi$-periodic for all $t>0$. Notice that $u_{\epsilon}(\cdot,t)$ is even with respect to $x=0$ and $x=\pi$ (by uniqueness, since so is $g_\epsilon$). We deduce that the derivative $v_{\epsilon}:=\partial_x u_{\epsilon}$ solves
$$
  \left\{\begin{array}{ll}
           \partial_t v_{\epsilon}-\partial_{xx} v_{\epsilon}=0 &\text{ in }(0,\pi)\times (0,+\infty),
           \smallskip \\
           v_{\epsilon}(0,t)=v_{\epsilon}(\pi,t)=0 &\text{ for all $t> 0$}.
         \end{array}\right.
$$
Moreover, by the assumptions on $g_{\epsilon}$, we have that
$v_{\epsilon}(x,0)=\partial_x u_{\epsilon} (x,0)=g_\epsilon'(x)\leq0$ for all $x\in [0,\pi]$.
Hence, the maximum principle yields
\begin{equation}
  \label{eq:v eps bigger zero}
   v^{\epsilon}\leq0\quad \text{ in }[0,\pi]\times[0,+\infty).
\end{equation}

We now claim that,  for every given $t>0$ and $m=0,1,2,\ldots$, 
\begin{equation}
 \label{eq:u eps conv to Gamma}
  u_{\epsilon}(\cdot,t)\text{ converges to } \Gamma(\cdot,t) \text{ in $C^m([0,\pi])$ as $\epsilon\downarrow 0$}.
\end{equation}
This follows by combining the uniform continuity in 
$[0,\pi]$ of $\Gamma(\cdot,t)$ and of all its derivatives with the fact that, for every $\delta>0$, we have
$$
  \big|\partial^m_x u_{\epsilon}(x,t)
  -\partial^m_x \Gamma(x,t)\big|
  \leq \int_{-\epsilon}^{\epsilon}dy\,
  \big|\partial^m_x\Gamma(x-y,t)
  -\partial^m_x\Gamma(x,t)\big|
  |g_{\epsilon}(y)|\leq \delta
$$
if $\epsilon$ is chosen small enough. 

Now, from \eqref{eq:v eps bigger zero} and \eqref{eq:u eps conv to Gamma} we conclude that
$$
  \partial_x\Gamma(x,t)=\lim_{\epsilon\downarrow 0}\partial_x u_{\epsilon}(x,t)=\lim_{\epsilon\downarrow 0} v_{\epsilon}(x,t)\leq 0
$$
for all $x\in [0,\pi]$ and $t>0$.

Finally, we show that $\Gamma(\cdot,t)$ is decreasing in $(0,\pi)$. For this, notice that we already proved that, for all $t>0$, $\Gamma(\cdot,t)$ is even with respect to $x=0$ and $x=\pi$. Therefore, given any  $t_0>0$, we see that $v:=\partial_x\Gamma$ solves
\begin{equation}\label{eq:fund property of v}
  \left\{\begin{array}{ll}
           \partial_t v-\partial_{xx}v=0 &\text{ in }(0,\pi)\times (t_0,+\infty),
           \smallskip \\
           v(0,t)=v(\pi,t)=0 & \text{ for all $t>t_0$},
           \smallskip \\
           v(\cdot,t_0)=\partial_x\Gamma(\cdot,t_0)\leq0 &\text{ in $[0,\pi]$}.
         \end{array}\right.
\end{equation}
Now, the strong maximum principle yields that either $v<0$ in $(0,\pi)\times(t_0,+\infty)$ or $v\equiv0$ in 
$[0,\pi]\times[t_0,+\infty)$. The proof is now complete, since $v\equiv0$ is absurd ---it would give that $\Gamma(\cdot,t)$ is constant in $x$, clearly contradicting the fact that $\Gamma$ is the fundamental solution of the heat equation with periodic boundary conditions, as shown in the beginning of this appendix.
\end{proof}

\section*{Acknowledgments}

We thank M. Ritor\'e, E. Valdinoci, and T. Weth for sharing some references and information on the topic of the paper.



\end{document}